\documentclass[a4paper,12pt]{article}

\usepackage{amsmath}
\usepackage{a4wide}
\usepackage{amsfonts}
\newenvironment{pf}{\textbf{Proof:}}{\hspace{\stretch{1}}$\square$}

\newtheorem{thm}{Theorem}

\newtheorem{df}{Definition}
\newtheorem{pr}{Proposition}

\newtheorem{dthm}{Definition-Theorem}
\title{Poisson and Diffusion Approximation of Stochastic Schr\"odinger Equations with Control}
\author{Cl\'ement PELLEGRINI\\
\vspace{-0,3cm}\scriptsize{Institut C.Jordan}\\
\vspace{-0,3cm}\scriptsize{Universit\'e C.Bernard, Lyon 1}\\
\vspace{-0,3cm}\scriptsize{21, av Claude Bernard}\\
\vspace{-0,3cm}\scriptsize{69622 Villeurbanne Cedex}\\
\vspace{-0,3cm}\scriptsize{France}\\
\vspace{-0,3cm}\scriptsize{e-mail: pelleg@math.univ-lyon1.fr}}

\begin{document}

\maketitle

\begin{abstract}
"Quantum trajectories" are solutions of stochastic differential
equations. Such equations are called ``Stochastic Schr\"odinger
Equations'' and describe random phenomena in continuous
measurement theory of Open Quantum System. Many recent
developments deal with the control theory in such model
(optimization, monitoring, engineering...). In this article,
stochastic models with control are mathematically and physically
justified as limit of concrete discrete procedures called
``Quantum Repeated Measurements''. In particular, this gives a
rigorous justification of the Poisson and diffusion approximation
in quantum measurement theory with control. Furthermore we
investigate some examples using control in quantum mechanics.

\end{abstract}

\section*{Introduction}

One of the topic in Quantum Open System theory concerns the study
of the evolution of a small quantum system $\mathcal{H}_0$
undergoing an indirect and continuous measurement (the small
system is in contact with environment and the measurement is
performed on the environment). In this context, the evolution of
the system is usually described by classical stochastic
differential equations called \textit{"Stochastic Schr\"odinger
Equations"}. Essentially, two kind of equations are considered:
\begin{enumerate}
 \item  The ``diffusive
equation'' (Homodyne detection experiment) is given by
\begin{equation}\label{DIF}d\rho_t=L(\rho_t)dt+[\rho_t C^\star+ C\rho_t-Tr\left(\rho_t(C+C^\star)
\right)\rho_t]dW_t \end{equation} where $W_t$ describes a
one-dimensional Brownian motion.

\item The ``jump equation'' (Resonance fluorescence experiment) is
\begin{equation}\label{POIS}d\rho_t=L(\rho_t)dt+\left[\frac{\mathcal{J}(\rho_t)}{Tr[\mathcal{J}(\rho_t)]}-\rho_t\right](d\tilde{N}_t-Tr[\mathcal{J}(\rho_t)]dt)
\end{equation}
where $\tilde{N}_t$ is a counting process with stochastic
intensity $\int_0^t Tr[\mathcal{J}(\rho_s)]ds$.
\end{enumerate}

Solutions of such equations are called \textit{"Quantum
Trajectory"}; they describe the evolution of the state of the
system perturbed by the continuous indirect measurement.
\smallskip

Recent progress and developments in quantum optics and quantum
information theory need a highest precision in experience process
using measurement \cite{MR2271425} (sensitivity, miniaturization,
optimization...). It imposed then quantum systems to be
controlled. Two types of controls are usually considered in
Quantum Mechanics: deterministic and stochastic control.

A laser monitoring a qubit, i.e a two level atom system, is a
basic example of deterministic control application. The action of
control is then characterized by the laser intensity (the term
deterministic is related to the fact that we consider the
intensity to be a deterministic function of time). Often, it is
called \textit{"Open Loop Control"}. Such experiences are used in
order to prepare systems in specific states for quantum computing.

The notion of stochastic control is, here, directly connected with
the procedure of measurement. Depending on the result of
measurement, a control operation is performed in order to modify
the evolution of the system. As a result of measurement is random
in Quantum Mechanics (one of the axiom of the theory), the
control becomes also random (it justifies the term of stochastic
control). This kind of control is particularly used in
engineering when some constraint of precision and optimization
must be followed. Usually such control is called \textit{"Closed
Loop Control"} or \textit{"Feedback Control"}.
\smallskip

From a theoretical point of view, an important question is to lay
out a mathematical setup to modelize the action of control. The
next step is to describe the evolution of controlled quantum
system.

Usually in the literature, in order to obtain and justify the
classical stochastic Schr\"odinger equations (\ref{DIF}) and
(\ref{POIS}), Quantum Filtering theory \cite{MR2042615} or
Instrumental Process theory \cite{MR2124562} are used. Such
techniques are based on the Hilbertian formalism of Quantum
Mechanics and on the theory of Stochastic Quantum Calculus. It
uses heavy analytic machinery and all the subtleties of the non
commutative character of quantum probability (conditionnal
expectation in Von Neumann Algebra, partially observed system...).
The starting point is the description of interaction between
system and environment in terms of quantum stochastic differential
equations (also called \textit{Hudson Parthasarathy Equations}
\cite{MR1164866}). In order to apply such theory in the control
setup, a theory in adequacy with the non commutative character
have to be introduced. Even if it is satisfied, the derivation
and the obtaining of stochastic Schr\"odinger equations with
control is far from being obvious and intuitive (see
\cite{bouten-2006}) and there are less rigorous results.

Recently, in the framework of the description of the interaction
of a small system with environment (without measurement), in
\cite{AP}, the authors have introduced a discrete model of
interactions: \textit{"Quantum Repeated Interactions"}. The basic
model is the one of a small system $\mathcal{H}_0$ in contact
with an infinite chain of quantum system
$\bigotimes_{j=1}^\infty\mathcal{H}$. Ones after the others each
copy of $\mathcal{H}$ interacts with $\mathcal{H}_0$ during a
time $h$.

Such approach of open quantum system yields a "good" and "useful"
approximation model of continuous-time interaction models. Indeed
by rescaling this interaction with respect to $h$, it is shown
that models of interaction (described by stochastic quantum
differential equations) can be obtained as continuous limits (h
goes to zero) of discrete models. In the measurement setup, this
approach has been adapted in \cite{Poisson} and \cite{diffusion}.
In these articles, it is then shown that classical quantum
trajectories (solutions of equations (\ref{DIF}) and
(\ref{POIS})) can be obtained as continuous limits of discrete
models of quantum measurement called \textit{"Quantum Repeated
Measurements"} (also called \textit{"Discrete Quantum Trajectory
theory"}). The idea of discrete indirect measurements consists in
performing a measurement of an observable of $\mathcal{H}$ after
each interaction between $\mathcal{H}_0$ and a copy of
$\mathcal{H}$. A discrete quantum trajectory is then a discrete
random process describing the evolution of the state of
$\mathcal{H}_0$ undergoing such repeated measurements. In this
case, the approach of the theory of stochastic Schr\"odinger
equations via approximation results is essentially based on
classical probability theory (there are no problem of
commutativity).

The main aim of this article is to adapt such technique in the
framework of control. In this article, we present the notion of
control in the model of quantum repeated interactions. In this
setup, quantum repeated measurements give rise to a discrete
models of indirect quantum measurement with control. By adapting
convergence results of \cite{Poisson} and \cite{diffusion}, we
obtain the description of stochastic Schr\"odinger equations with
control. With this approach, all the problem concerning non
commutativity are avoided
 and the physical justification of stochastic models is rigorous and intuitive.\\

This article is structured as follows.

The first section is devoted to present discrete models of quantum
measurement with control theory. We remember the mathematical
model of quantum repeated interactions. Next, we introduce an
appropriate notion of control in this setup and by introducing
the measurement principle, we obtain the description of
\textit{discrete quantum trajectory with control}. Next, in order
to prepare final convergence results, we adapt and enlarge the
asymptotic assumptions presented in \cite{AP} to the context of
control. To investigate such problems, we focus on a central case
in physical application: a two-level atom in contact with a spin
chain.

The second section is then devoted to continuous models. By
applying the asymptotic assumptions on the two-level atom model of
Section $1$, we obtain two different discrete evolution equations
(in asymptotic form) describing the evolution of the state of
$\mathcal{H}_0$. Each evolution equation describes the evolution
of a discrete quantum trajectory with control for a specific
observable. For each equation, we investigate the continuous
limit equation and we show the convergence.

In the last section, we present some applications of continuous
models. In a first time, we investigate a model of a deterministic
control: an atom monitored by a laser. By modeling a suitable
interaction discrete model and by adapting the result of Section
$2$, we obtain a continuous stochastic model for this concrete
example. In a second time we present a use of stochastic control:
"Optimal Control theory". We adapt classical results concerning
this theory in the quantum language.

\section{Discrete Controlled Quantum Trajectories}

This section is devoted to the presentation of the model of
discrete quantum trajectories in presence of external control.
Here, we present a natural way for modeling control theory.

\subsection{Repeated Quantum Measurements with Control}

 In a first time,
let us remember the general context of quantum repeated
interactions.

 A small system, represented by
a Hilbert space $\mathcal{H}_0$, is in contact with an infinite
chain of identical independent quantum systems. Each piece of the
environment is represented by $\mathcal{H}$ and interacts, one
after the others, with $\mathcal{H}_0$ during a time interval of
length $h$. For example, a copy of $\mathcal{H}$ can represent an
incoming photon or a measurement apparatus...

The space describing the first interaction between $\mathcal{H}_0$
and $\mathcal{H}$ is defined by the tensor product
$\mathcal{H}_0\otimes\mathcal{H}$. The evolution is given by a
self-adjoint operator $H_{tot}$ on the tensor product. This
operator is called the total Hamiltonian and its general form is
$$H_{tot}=H_0\otimes I+I\otimes H+H_{int}$$
where the operators $H_0$ and $H$ are the free Hamiltonians of
each system. The operator $H_{int}$ represents the Hamiltonian of
interaction. This allows to define a unitary-operator
$$U=e^{ih\,H_{tot}},$$
and the evolution of states of $\mathcal{H}_0\otimes\mathcal{H}$,
in the Schr\"odinger picture, is given by
$$\rho\mapsto U\,\rho\, U^\star.$$
 After this first interaction, a second copy of $\mathcal{H}$ interacts with
$\mathcal{H}_0$ in the same fashion and so on.

As the chain is supposed to be infinite, the whole sequence of
interactions is described by the state space:
\begin{equation}
\mathbf{\Gamma}=\mathcal{H}_0\otimes\bigotimes_{k\geq
1}\mathcal{H}_k
\end{equation}
where $\mathcal{H}_k$ denotes the $k$-th copy of $\mathcal{H}$.
The countable tensor product $\bigotimes_{k\geq 1}\mathcal{H}_k$
means the following. Consider that $\mathcal{H}$ is of finite
dimension and that $\{X_0,X_1,\ldots,X_n\}$ is a fixed orthonormal
basis of $\mathcal{H}$. The orthogonal projector on
$\mathbb{C}X_0$ is denoted by $\vert X_0\rangle\langle X_0\vert$.
This is the ground state (or vacuum state) of $\mathcal{H}$. The
tensor product is taken with respect to $X_0$ (for more details,
see \cite{AP}).

The unitary evolution describing the k-th interaction is given by
$U_k$ which acts as $U$ on $\mathcal{H}_0\otimes\mathcal{H}_k$,
whereas it acts like the identity operator on the other copies.
If $\rho$ is a state on $\mathbf{\Gamma}$, the effect of the
$k$-th interaction is then:
$$\rho\mapsto U_k\,\rho\, U_k^\star$$
 Hence the result of the $k$ first interactions is described by the operator $V_k$ on
$\mathcal{B}(\mathbf{\Gamma})$ defined by the recursive formula:
 \begin{equation}\left\{\begin{array}{ccc} V_{k+1}&=&U_{k+1}V_k \\
  V_0&=&I\end{array}\right.\end{equation}
  and the evolution of states is given by
 $$\rho\mapsto V_k\,\rho\, V_k^\star.$$

 In this context, a main feature of this article is to present measurement and control theory.
  Let start by describing the control theory. An action of control consist
  in modifying the interaction at
  each new step depending on the previous step
  (this condition allows further to introduce stochastic control). Therefore
if $U_k$ is the unitary-operator describing the $k$-th
interaction, it depends then on the time of interaction and
 on a parameter $u_{k-1}$ which gives account of the control.
 Likewise this parameter depends on the interaction time; the operator $U_k$ is then denoted by
 $U_k(h,u_{k-1}(h))$.

  The whole sequence $\mathbf{u}=(u_k(h))$ is
 called the "control strategy". In term of $\mathbf{u}$, the $k$ first interactions are
 then
 described by the unitary-operator $V_k^\mathbf{u}$:
\begin{equation}
 V_k^\mathbf{u}=U_k(h,u_{k-1}(h))\,U_{k-1}(h,u_{k-2}(h))\,\ldots \,U_1(h,u_0(h)).
\end{equation}
Finally, the evolution in presence of control is given by
\begin{equation}\label{evol}
 \rho\mapsto V_k^\mathbf{u}\,\,\rho\,\, V_k^{\mathbf{u}\star}.
\end{equation}
Before to give a complete definition of control strategies,
 we have to introduce the repeated measurement model.
 \smallskip

Let us describe the basic procedure on each system of the chain.
Let $A$ be any observable on $\mathcal{H}_k$ with spectral
decomposition $A=\sum_{j=0}^{p}\lambda_j P_j$, consider its
natural ampliation as an observable on $\mathbf{\Gamma}$ by:
\begin{equation}\label{ampliation}A^k:=\bigotimes_{j=0}^{k-1}I\otimes A\otimes\bigotimes_{j\geq k+1}I\end{equation}
 The accessible data are the eigenvalues of $A^k$ and the result of the observation is random.
 If $\rho$ is any state on $\mathbf{\Gamma}$, we observe $\lambda_j$ with probability
 $$P[\textrm{to observe}\,\,\lambda_j]=Tr[\,\rho\,P^k_j\,],\,\,\,\,j\in\{0,\ldots,p\},$$
 where the operator $P^k_j$ corresponds to the ampliation of the eigenprojector $P_j$ in the same way as $(\ref{ampliation})$. If we have observed the eigenvalue $\lambda_j$, the ``projection'' postulate called
``wave packet reduction'' imposes the state after the measurement
to be
 $$\rho_j=\frac{P^k_j\,\rho\,P^k_j }{Tr[\,\rho\,P^k_j\,]}.$$

 $\textbf{Remark}$: This corresponds to the new reference state of our system.
 Another measurement of the same observable $A^k$ (with respect to this state)
 should give $P[\textrm{to observe}\,\,\lambda_j]=1$.
 Hence only one measurement give a significant information;
 it justifies a principle of repeated interactions.\\

 Quantum repeated measurements with control are the combination
 of this previous principle and the successive interactions $(\ref{evol})$.
 After each interaction, a quantum measurement induces a random modification of the state of the system. It defines then a discrete process which is called
 ``discrete controlled quantum trajectory''. The description is as follows.

 The initial state on $\mathbf{\Gamma}$ is chosen to be
 $$\mu=\rho\otimes\bigotimes_{j\geq 1}\beta_j$$
 where $\rho$ is any state on $\mathcal{H}_0$ and each $\beta_i=\beta$ where $\beta$ is
  any state on $\mathcal{H}$. The state after $k$ interactions
  is denoted by $\mu_k^\mathbf{u}$, we have: $$\mu_k^\mathbf{u}=V_k^\mathbf{u}\,\,\mu\,\,
V^{\mathbf{u}\star}_k.$$

 The probability space describing the experience
 is $\Sigma^{\mathbb{N}^\star}$ where $\Sigma=\{0,\ldots,p\}$.
 The integers $i$ correspond to the indexes of the eigenvalues of $A$.
We endow $\Sigma^{\mathbb{N}^\star}$ with the cylinder
$\sigma$-algebra generated by the cylinder sets:
 $$\Lambda_{i_1,\ldots,i_k}=\{\omega\in\Omega^\mathbb{N}/\omega_1=i_1,\ldots,\omega_k=i_k
\}.$$

 Remarking that for all $j$, the unitary operator $U_j$ commutes with all $P^k$
for all $k<j$. For any set $\{i_1,\ldots,i_k\}$, we can define
the following operator:
\begin{eqnarray*}\tilde{\mu}^\mathbf{u}_k(i_1,\ldots,i_k)&=&I\otimes P_{i_1}\otimes
\ldots\otimes P_{i_k}\otimes I\ldots
\,\,\mu_k^\mathbf{u}\,\,I\otimes P_{i_1}\otimes \ldots\otimes
P_{i_k}\otimes
I\ldots\\
&=&P^k_{i_k}\ldots
P^1_{i_1}\,\,\,\mu_k^\mathbf{u}\,\,\,P^1_{i_1}\ldots P^k_{i_k}.
\end{eqnarray*}
This is the non-normalized state corresponding to the successive
observation of $\lambda_{i1},\ldots,\lambda_{ik}$. The probability
to observe these eigenvalues is
$$P[\textrm{to\,\,observe}\,\,\lambda_{i1},\ldots,\lambda_{ik}]=Tr[\tilde{\mu}_k^\mathbf{u}(i_1,\ldots,i_k)].$$
By putting
$$P[\Lambda_{i_1,\ldots,i_k}]=P[\textrm{to\,\,observe}\,\,\lambda_{i1},\ldots,\lambda_{ik}],$$
 it defines a probability measure on the cylinder sets of $\Sigma^{\mathbf{N}^\star}$ which satisfies the Kolmogorov Consistency Criterion. It defines then a unique probability measure on $\Sigma^{\mathbf{N}^\star}$. The discrete quantum trajectory with control strategy $\mathbf{u}$ on $\mathbf{\Gamma}$ is described by the following random sequence of states:
 $$\begin{array}{cccc}\tilde{\rho}_k^\mathbf{u}: & \Sigma^{\mathbf{N}^\star} &\longrightarrow&
\mathcal{B}(\mathbf{\Gamma})\\
 & \omega & \longmapsto & \tilde{\rho}_k^\mathbf{u}(\omega_1,\ldots,
\omega_k)=\frac{\tilde{\mu}^\mathbf{u}_k(\omega_1,\ldots,
\omega_k)}{Tr[\tilde{\mu}^\mathbf{u}_k(\omega_1,\ldots,
\omega_k))]}.
\end{array}$$
From this description, the following result is obvious.
\begin{pr}\label{p}
Let $\mathbf{u}$ be any strategy and
$(\tilde{\rho}^\mathbf{u}_k)$ be the above random sequence of
states we have for all $\omega\in\Sigma^{\mathbb{N}}$:
$$\tilde{\rho}^\mathbf{u}_{k+1}(\omega)=\frac{P^{k+1}_{\omega_{k+1}}\,
U_{k+1}(h,u_{k}(h))\,\,\,\,\tilde{\rho}_k^\mathbf{u}(\omega)\,\,\,\,U^\star_{k+1}(h,u_{k}(h))\,
P^{k+1}_{\omega_{k+1}}}{Tr\left[\,\tilde{\rho}_k^\mathbf{u}(\omega)\,\,U^\star_{k+1}(h,u_{k}(h))\,
P^{k+1}_{\omega_{k+1}}\,U_{k+1}(h,u_{k}(h))\right]}.$$
\end{pr}

Now at this stage, we can make precise the definition of control
strategies which correspond to the case of deterministic or
stochastic control mentioned in the Introduction.

\begin{df}
 Let $\mathbf{u}=(u_k(h))$  be a control strategy and let $(\tilde{\rho}^{\mathbf{u}}_k)$ be a quantum trajectory.
\begin{enumerate}
 \item {If there exists some function $u$ from $\mathbb{R}$ to $\mathbb{R}^{n}$ such that for all $k:$ $$u_{k}(h)=u(kh),$$
  the control strategy is called deterministic. It is also called ``open loop control''.}
\item{If there exists some function $u$ from $\mathbb{R}\times\mathcal{B}(\mathbf{\Gamma})$
 to $\mathbb{R}^{n}$ such that for all $k:$ $$u_{k}(h)=u(kh,\tilde{\rho}^{\mathbf{u}}_k),$$ the control strategy is called Markovian.
  It is also called ``closed loop control'' or ``feedback control''. If for all $k$ we have $u_{k}(h)=u(\tilde{\rho}^\mathbf{u}_k)$, this is an homogeneous Markovian strategy.}
\end{enumerate}
\end{df}

 The following theorem is an easy
consequence of Proposition $\ref{p}$ and of the previous
Definition.

\begin{thm}
 For all control strategy $\mathbf{u}$, the sequence $(\tilde{\rho}_n^\mathbf{u})_n$ is a
 non homogeneous Markov chain valued on the
set of states
  of $\mathbf{\Gamma}$. It is described as follows:
\begin{eqnarray*}
P\left[\tilde{\rho}_{n+1}^\mathbf{u}=\mu/\tilde{\rho}_n^\mathbf{u}=\theta_n,\ldots,\tilde{\rho}^\mathbf{u}_0=\theta_0\right]
=P\left[\tilde{\rho}_{n+1}^\mathbf{u}=\mu/\tilde{\rho}_n^\mathbf{u}=\theta_n\right].
\end{eqnarray*}If $\tilde{\rho}_n^\mathbf{u}=\theta_n$ then $\tilde{\rho}_{n+1}^\mathbf{u}$ takes one
of the values:
$$ \mathcal{H}_i^{\mathbf{u},n+1}(\theta_n)=\frac{
P_i^{n+1}\,U_{n+1}(h,u_{n}(h))\,\,\theta_n\,\,U^\star_{n+1}(h,u_{n}(h))\,
P_i^{n+1}}{Tr\left[\,U_{n+1}(h,u_{n}(h))\,\,\theta_n\,\,U^\star_{n+1}(h,u_{n}(h))\,
P_i^{n+1}\,\right]},\,\,\,\,i=0,\ldots,p,$$ with probability
$Tr\left[\,U_{n+1}(h,u_{n}(h))\,\theta_n\,U^\star_{n+1}(h,u_{n}(h))\,P_i^{n+1}\,\right].$

The discrete process $(\rho_k^\mathbf{u})$ is called a controlled
Markov chain.
\end{thm}

\begin{pf}
 Property of being a Markov chain comes from the fact that a control strategy is either deterministic or Markovian.
  For the two cases, the conclusion is obvious from the description of Proposition $\ref{p}$.
\end{pf}\\

In general, one is only interested in the reduced state of the
small system. This state is given by the partial trace operation.
Let us recall what partial trace is. Let $\mathcal{Z}$ be any
Hilbert space, the notation $Tr_{\mathcal{Z}}[W]$ corresponds to
the trace of any trace-class operator $W$ on $\mathcal{Z}$.

\begin{dthm} Let $\mathcal{H}$ and $\mathcal{K}$ be any Hilbert spaces.
Let $\alpha$ be a state on the tensor product
$\mathcal{H}\otimes\mathcal{K}$. There exists a unique state
$\eta$ on $\mathcal{H}$ which is characterized by the property:
 $$ Tr_{\mathcal{H}}[\,\eta
\,X\,]=Tr_{\mathcal{H}\otimes\mathcal{K}}[\,\alpha (X\otimes
I)\,].$$ for all $X\in\mathcal{B}(\mathcal{H})$. The state $\eta$
is called the partial trace of $\alpha$ on $\mathcal{H}$ with
respect to $\mathcal{K}$.
\end{dthm}

For any state $\alpha$ on $\mathbf{\Gamma}$, denote
$\mathbf{E}_0[\alpha]$ the partial trace of $\alpha$ on
$\mathcal{H}_0$ with respect to
$\bigotimes_{k\geq1}\mathcal{H}_k$. We then define the discrete
controlled quantum trajectory on $\mathcal{H}_0$ as follows. For
all $\omega\in\Sigma^{\mathbb{N}^\star}$:
\begin{equation}\label{TrPart}\rho_n^\mathbf{u}(\omega)=\mathbf{E}_0[\tilde{\rho}_n^\mathbf{u}(\omega)].\end{equation}
\textbf{Remark}: We adapt Definition $1$ by considering Markovian strategy defined on $\mathbb{R}\times\mathcal{B}(\mathcal{H}_0)$.\\
An immediate consequence of Theorem $1$ is the following result.

\begin{thm}\label{RaS}
For all control strategy $\mathbf{u}$, the random sequence defined
by formula $(\ref{TrPart})$ is a non-homogeneous controlled
Markov chain with values in the set of states on $\mathcal{H}_0$.
If $\rho^\mathbf{u}_n=\chi_n$ then $\rho^\mathbf{u}_{n+1}$ takes
one of the values:
$$\mathbf{E}_0\left[\frac{I\otimes P_i\,\,\tilde{U}_{n+1}(h,u_{n}(h))\,(\chi_n\otimes\beta)\,\tilde{U}^\star_{n+1}(h,u_{n}(h))\,\,I\otimes
P_i}{Tr[\,\tilde{U}_{n+1}(h,u_{n}(h))\,(\chi_n\otimes\beta)\,\tilde{U}^{\star}_{n+1}(h,u_{n}(h))\,I\otimes
P_i]}\right] \,\,i=0\ldots p$$ with probability
$Tr\left[\tilde{U}_{n+1}(h,u_{n}(h))(\chi_n\otimes\beta)\tilde{U}^\star_{n+1}(h,u_{n}(h))\,P_i\right]$.
\end{thm}
\textbf{Remark}: Let us stress that: $$\frac{(I\otimes
P_i)\,U\,(\chi_n\otimes\beta)\,U^\star\,(I\otimes
P_i)}{Tr[\,U\,(\chi_n\otimes\beta)\,U^{\star}\,(I\otimes P_i)]}$$
is a state on $\mathcal{H}_0\otimes\mathcal{H}$. In this
situation, the notation $\mathbf{E}_0$ denotes the partial trace
on $\mathcal{H}_0$ with respect to $\mathcal{H}$. Moreover for
each n, the operator $\tilde{U}_n$, which appears in the
description of the transition of the Markov chain
$(\rho^\mathbf{u}_n)$, acts on $\mathcal{H}_0\otimes\mathcal{H}$
as the operators $U_n$ on
$\mathcal{H}_0\otimes\mathcal{H}_n$\\

With the description of Theorem 2, we can express a discrete
evolution equation describing the discrete quantum trajectory
$(\rho_k^\mathbf{u})$. By putting
$$\mathcal{L}^{\mathbf{u},k}_i(\rho)=\mathbf{E}_0\left[I\otimes P_i\,\tilde{U}_{k}(h,u_{k-1}(h))\,(\rho\otimes\beta)\,\tilde{U}^\star_{k}(h,u_{k-1}(h))\,I\otimes
P_i\right] \,\,i=0\ldots p,$$ and
$\mathbf{1}_i^{k}(\omega)=\mathbf{1}_i(\omega_{k})$ for all
$\omega\in\Sigma^{\mathbb{N}^\star}$, the discrete process
$(\rho_k^\mathbf{u})$ then satisfies
\begin{equation}\label{discrete}
\rho_{k+1}^\mathbf{u}(\omega)=
\sum_{i=0}^p\frac{\mathcal{L}^{k+1}_i(\rho^\mathbf{u}_k(\omega))}{Tr[\mathcal{L}^{k+1}_i(\rho^\mathbf{u}_k(\omega))]}\mathbf{1}_i^{k+1}(\omega)
\end{equation}
for all $\omega\in\Sigma^\mathbb{N}$ and all $k>0$.

The following section is devoted to the deeply study of the
equation $(\ref{discrete})$ in a particular case of a two-level
system in interaction with a spin chain. Next we come into the
question of asymptotic assumptions.

\subsection{A Two-Level Atom}

 The physical situation is described by
 $\mathcal{H}_0=\mathcal{H}=\mathbb{C}^2$. In this case, an
 observable $A$ has two different eigenvalues: $A=\lambda_0P_0+\lambda_1P_1$ (the case which only
 one eigenvalue is not interesting). The equation $(\ref{discrete})$
 is reduced to:
\begin{equation}
\rho^\mathbf{u}_{k+1}(\omega)=\frac{\mathcal{L}^{\mathbf{u},k+1}_{0}(\rho_k^\mathbf{u}(\omega))}{p^\mathbf{u}_{k+1}}\mathbf{1}_0^{k+1}(\omega)+\frac{\mathcal{L}^{\mathbf{u},k+1}_{1}
(\rho_k^\mathbf{u}(\omega))}{q^\mathbf{u}_{k+1}}\mathbf{1}_1^{k+1}(\omega).
\end{equation}
where
$p^\mathbf{u}_{k+1}=Tr[\mathcal{L}^{\mathbf{u},k+1}_{0}(\rho_k^\mathbf{u})]=1-q^\mathbf{u}_{k+1}$.
Let now introduce the centered and normalized random variable
$$X_{k+1}=\frac{\mathbf{1}_1^{k+1}(\omega)-q^\mathbf{u}_{k+1}}{\sqrt{q^\mathbf{u}_{k+1}p^\mathbf{u}_{k+1}}}.$$
We define the associated filtration on $\{0,1\}^{\mathbf{N}}$:
$$\mathcal{F}_k=\sigma(X_i,i\leq k).$$ So by construction we have
$\mathbf{E}[X_{k+1}/\mathcal{F}_k]=0$ and
$\mathbf{E}[X_{k+1}^2/\mathcal{F}_k]=1$. In terms of $(X_k)$ the
discrete controlled quantum trajectory satisfies:
\begin{eqnarray}\label{discreteequation}
\rho_{k+1}^\mathbf{u}&=&\mathcal{L}^{\mathbf{u},k+1}_{0}(\rho_k^\mathbf{u})+\mathcal{L}^{\mathbf{u},k+1}_{1}(\rho_k^\mathbf{u})\nonumber\\&&+\left[-\sqrt{\frac{q^\mathbf{u}_{k+1}}{p^\mathbf{u}_{k+1}}}\mathcal{L}^{\mathbf{u},k+1}_{0}(\rho_k^\mathbf{u})
+\sqrt{\frac{p^\mathbf{u}_{k+1}}{q^\mathbf{u}_{k+1}}}\mathcal{L}^{k+1}_{1}(\rho_k^\mathbf{u})\right]X_{k+1}.
\end{eqnarray}

To give more sense to the equation (\ref{discreteequation}), we
have to express the terms
$\mathcal{L}^{\mathbf{u},k+1}_{i}(\rho_k^\mathbf{u})$ in a more
explicit way. For this, we introduce a specific basis. Let
$(X_0=\Omega,X_1=X)$ be an orthonormal basis of
$\mathcal{H}_0=\mathcal{H}=\mathbb{C}^2$. For the space
$\mathcal{H}_{0}\otimes\mathcal{H}$,
  we consider the following basis
 $$\Omega\otimes\Omega,X\otimes\Omega,\Omega\otimes
X,X\otimes X.$$ In this basis, the unitary operator can be
written by blocks as a $2\times2$ matrix:
$$U_{k+1}(h,u_{k}(h))=\left( \begin{array}{cc}  L_{00}(kh,u_{k}(h)) & L_{01}(kh,u_{k}(h)) \\
  L_{10}(kh,u_{k}(h)) & L_{11}(kh,u_{k}(h))
 \end{array}\right) $$
 where each $L_{ij}(kh,u_{k}(h))$ are operators on $\mathcal{H}_0$. The reference state $\beta$ of $\mathcal{H}$ is:
$$\beta=\vert\Omega\rangle\langle\Omega\vert.$$
The terms $\mathcal{L}_i^{\mathbf{u},k+1}(\rho_k^\mathbf{u})$
depends also on the expression of the eigenprojectors of the
observable $A$. If the eigenprojector $P_i$ is expressed as
$P_i=\left(\begin{array}{cc}p_{00}^i&p_{01}^i\\p_{10}^i&p_{11}^i\end{array}\right)$
in the basis $(\Omega,X)$ of $\mathcal{H}$, we have:
\begin{eqnarray}\label{Lij}
 \mathcal{L}_i^{\mathbf{u},k+1}(\rho_k^\mathbf{u})&=&p_{00}^iL_{00}(kh,u_{k}(h))
 \,\rho_k^\mathbf{u}
\,L_{00}^\star(kh,u_{k}(h))+p_{01}^iL_{00}(kh,u_{k}(h))\,\rho_k^\mathbf{u}
\,L_{10}^\star(kh,u_{k}(h))\nonumber\\
&&+\,p_{10}^iL_{10}(kh,u_{k}(h))\,\rho_k^\mathbf{u}
\,L_{00}^\star(kh,u_{k}(h))+p_{11}^iL_{10}(kh,u_{k}(h))\,\rho_k^\mathbf{u}
\,L_{10}^\star(kh,u_{k}(h))\nonumber\\
\end{eqnarray}

As the unitary evolution depends on the time length interaction
$h$, the discrete quantum trajectory $(\rho_k^\mathbf{u})$
depends on $h$. In Section $2$, this dependence allows us to
consider continuous time limit ($h\rightarrow0$) of the discrete
processes $(\rho_k^\mathbf{u})$. The next section is devoted to
present the asymptotic ingredients necessary to obtain such
convergence results.

\subsection{Description of Asymptotic}

In this section, we present suitable asymptotic for the
coefficients of the unitary operators $U_{k}(h,u_k(h))$ in order
to have an effective continuous time limit from discrete quantum
trajectories. Let $h=1/n$ be the length time of interaction, we
have for $(U_k)$
$$U_{k+1}(n,u_{k}(n))=\left( \begin{array}{cc}  L_{00}(k/n,u_{k}(n)) & L_{01}(k/n,u_{k}(n)) \\
  L_{10}(k/n,u_{k}(n)) & L_{11}(k/n,u_{k}(n))
 \end{array}\right),$$
In our context, the choice of the coefficients $L_{ij}$ is an
adaptation of the works of Attal-Pautrat in \cite{AP}. In their
work, they consider only evolution of the type
$$U_{k+1}(n)=\left( \begin{array}{cc}  L_{00}(n) & L_{01}(n) \\
  L_{10}(n) & L_{11}(n)
 \end{array}\right),$$
that is, homogeneous evolution without control. They have shown
that
 $$V_{[nt]}=U_{[nt]}(n)\ldots U_1(n)$$ converges (in operator algebra) to a
  non-trivial process $V_t$ (solution of a quantum stochastic differential equation),
   only if the coefficients $L_{ij}(n)$ obey certain normalization. In their case,
   these coefficients must be of the form
\begin{eqnarray}\label{AAPP}L_{00}(n)&=&I+\frac{1}{n}\left(-iH_0-\frac{1}{2}CC^\star\right)+\circ\left(\frac{1}{n}\right)\\
  L_{10}(n)&=&\frac{1}{\sqrt{n}}C+\circ\left(\frac{1}{n}\right),
  \end{eqnarray}
where $H_0$ is the Hamiltonian of $\mathcal{H}_0$ and $C$ is any
operator on $\mathbb{C}^2$. Hence, in the control context, the
coefficients $L_{ij}(k/n,u_k(n))$ must follow similar
expressions. Let $k$ be fixed, we put
\begin{eqnarray}
 L_{00}(k/n,u_{k}(n))&=&I+\frac{1}{n}\left(-iH_k(n,u_k(n))-\frac{1}{2}C_k(n,u_k(n))C_k^\star(n,u_k(n))\right)+\circ\left(\frac{1}{n}\right)\nonumber\\
L_{00}(k/n,u_{k}(n))&=&\frac{1}{\sqrt{n}}C_k(n,u_k(n))+\circ\left(\frac{1}{n}\right)
\end{eqnarray}
where $H_k(n,u_k(n))$ is a self-adjoint operator and
$C_k(n,u_k(n))$ is an operator on $\mathbb{C}^2$. It is
straightforward that the expression $(\ref{AAPP})$ of
Attal-Pautrat is a particular case of the previous expression.
Finally, we suppose that there exist some function $H$ and $C$
such that
$$\begin{array}{ccccccccc}
H:&\mathbb{R}^+\times\mathbb{R}&\longrightarrow&\mathbb{H}_2(\mathbb{C})&\textrm{and}&C:&\mathbb{R}^+\times\mathbb{R}&\longrightarrow&\mathbb{M}_2(\mathbb{C})\\
&(t,s)&\longmapsto&H(t,s)&&&(t,s)&\longmapsto&C(t,s)
\end{array}$$
where $\mathbb{H}_2(\mathbb{C})$ designs the set of self-adjoint
operators on $\mathbb{C}^2$ and
\begin{eqnarray}
 H_k(n,u_k(n))&=&H(k/n,u_k(n))\nonumber\\
C_k(n,u_k(n)&=&C(k/n,u_k(n))
\end{eqnarray}
Furthermore we suppose that all the $\circ$ are uniform in $k$.

Now, we shall express the equation $(\ref{discreteequation})$ and
$(\ref{Lij})$ with these asymptotic assumptions. As it was
announced, we obtain two different behaviours depending of the
choice of the observable.

\begin{enumerate}
\item If the observable $A$ is diagonal in the basis $(\Omega,X)$, that is, it is of the form\\ $A=\lambda_0\left( \begin{array}{cc}  1 & 0 \\
  0 & 0
 \end{array}\right)+\lambda_1\left( \begin{array}{cc}  0 & 0 \\
  0 & 1
 \end{array}\right)$,
 we obtain the asymptotic for the probabilities
 \begin{eqnarray*}
p^\mathbf{u}_{k+1}(n)&=&1-\frac{1}{n}Tr\Big[\mathcal{J}(k/n,u_{k}(n))(\rho^\mathbf{u}_k(n))\Big]+\circ\left(\frac{1}{n}\right)\\
q^\mathbf{u}_{k+1}(n)&=&\frac{1}{n}Tr\Big[\mathcal{J}(k/n,u_{k}(n))(\rho^\mathbf{u}_k(n))\Big]+\circ\left(\frac{1}{n}\right)
\end{eqnarray*}
The discrete equation $(\ref{discreteequation})$ becomes
\begin{eqnarray}\label{poiss}
&&\hspace{-1,5cm}\rho^\mathbf{u}_{k+1}(n)-\rho^\mathbf{u}_k(n)\nonumber\\\hspace{-0,cm}&&\hspace{-0,9cm}=\frac{1}{n}L(k/n,u_{k}(n))(\rho^\mathbf{u}_k(n))+\circ(\frac{1}{n})\nonumber\\
&&\hspace{-0,5cm}+\left[\frac{\mathcal{J}(k/n,u_{k}(n))(\rho^\mathbf{u}_k(n))}{Tr\Big[\mathcal{J}(k/n,u_{k}(n))(\rho^\mathbf{u}_k(n)))\Big]}-\rho^\mathbf{u}_k(n)+\circ(1)\right]
\sqrt{q^\mathbf{u}_{k+1}(n)p^\mathbf{u}_{k+1}(n)} \,X_{k+1}(n)
\end{eqnarray}
where for all state $\rho$, we have defined
\begin{eqnarray}\label{LINDBLAD}
\mathcal{J}(t,s)(\rho)&=&C(t,s)\,\,\rho\,\,C^\star(t,s)\,\,\,\textrm{and}\nonumber\\
L(t,s)(\rho)&=&-i[H(t,s),\rho]-\frac{1}{2}\{C(t,s)C^\star(t,s),\rho\}+\mathcal{J}(t,s)(\rho).
\end{eqnarray}

\item If the observable $A$ is non diagonal in the basis $(\Omega,X)$, and if the eigenprojectors are express as $P_0=\left(
\begin{array}{cc}
p_{00} & p_{01} \\
  p_{10} & p_{11}
 \end{array}\right)$ and $P_1=\left( \begin{array}{cc}  q_{00} & q_{01} \\
  q_{10} & q_{11}
 \end{array}\right)$ we have
 \begin{eqnarray*}
p^\mathbf{u}_{k+1}&=&p_{00}+\frac{1}{\sqrt{n}}Tr\Big
[\rho^\mathbf{u}_k\,\,\big(p_{01}C(k/n,u_{k+1}(n))+p_{10}C^\star(k/n,u_{k}(n))\big)\Big]\nonumber\\&&
+\frac{1}{n}Tr\Big[\rho^\mathbf{u}_k\,p_{00}\,\,\big(C(k/n,u_{k}(n))+C^\star(k/n,u_{k}(n))\big)
\Big]
+\circ\left(\frac{1}{n}\right)\\
q^\mathbf{u}_{k+1}&=&q_{00}+\frac{1}{\sqrt{n}}Tr\Big[
\rho^\mathbf{u}_k\,\,\big(q_{01}C(k/n,u_{k}(n))+q_{10}C^\star(k/n,u_{k}(n))\big)\Big]\nonumber\\&&
+\frac{1}{n}Tr\Big[\rho^\mathbf{u}_k\,q_{00}\,\,\big(C(k/n,u_{k}(n))+C^\star(k/n,u_{k}(n))\big)\Big]+\circ\left(\frac{1}{n}\right).
\end{eqnarray*}
The discrete equation $(\ref{discreteequation})$ becomes
\begin{eqnarray}\label{diffF}
&&\hspace{-1,0cm}\rho^\mathbf{u}_{k+1}-\rho^\mathbf{u}_k=\nonumber\\&&
\frac{1}{n}L(k/n,u_{k}(n))(\rho^\mathbf{u}_k)+\circ\left(\frac{1}{n}\right)+
\Big[e^{i\theta}C(k/n,u_{k}(n))\rho^\mathbf{u}_k
+e^{-i\theta}\rho^\mathbf{u}_kC^\star(k/n,u_{k}(n))\nonumber\\&&
 -Tr\big[\rho^\mathbf{u}_k\,\,\big(e^{i\theta}C(k/n,u_{k}(n))+e^{-i\theta}C^\star(k/n,u_{k}(n))\big)\big]\,\rho^\mathbf{u}_k+\circ
(1)\Big]\frac{1}{\sqrt{n}}X_{k+1}(n)
\end{eqnarray}
where $\theta$ is a real parameter. This parameter can be
explicitly expressed with the coefficients of the eigenprojectors
$(P_i)$. By putting
$C_\theta(k/n,u_{k}(n))=e^{i\theta}C(k/n,u_{k}(n))$ we have the
same form for the equation $(\ref{diffF})$ for all $\theta$, then
we consider in the following that $\theta=0$. The expression of
$L$ is the same as $(\ref{LINDBLAD})$.
\end{enumerate}

In order to prepare the final convergence result, in each case, we
can define a process $(\rho_{[nt]})$ which satisfies
\begin{eqnarray}\label{DIScrete}
\rho_{[nt]}^\mathbf{u}&=&\rho_0+\sum_{i=0}^{[nt]-1}[\rho_{i+1}^\mathbf{u}-\rho_i^\mathbf{u}]\nonumber\\&=&\rho_0+\sum_{i=0}^{[nt]-1}
[\mathcal{L}^{\mathbf{u},i+1}_{0}(\rho_i^\mathbf{u})+\mathcal{L}^{\mathbf{u},i+1}_{1}(\rho_i^\mathbf{u})-\rho_i^\mathbf{u}]\nonumber\\
&&+\sum_{i=0}^{[nt]-1}\left[
-\sqrt{\frac{q^\mathbf{u}_{i+1}}{p^\mathbf{u}_{i+1}}}\,\mathcal{L}^{\mathbf{u},i+1}_{0}(\rho_i^\mathbf{u})+\sqrt{\frac{p^\mathbf{u}_{i+1}}{q^\mathbf{u}_{i+1}}}\,\mathcal{L}^{\mathbf{u},i+1}_{1}(\rho_i^\mathbf{u})\right]
X_{i+1}\nonumber\\
&=&\rho_0+\sum_{i=0}^{[nt]-1}\frac{1}{n}\mathcal{Y}(i/n,u_i(n),\rho^\mathbf{u}_i)+\sum_{i=0}^{[nt]-1}\mathcal{Z}(i/n,u_i(n),\rho^\mathbf{u}_i)X_{i+1}
\end{eqnarray}
for some functions $\mathcal{Y}$ and $\mathcal{Z}$ which depend
on the description $(\ref{poiss})$ or $(\ref{diffF})$.

Depending on the choice of observable, in the next section, we
show that the process $(\rho_{[nt]})$ converges to a solution of a
particular stochastic differential equation.

\section{Convergence to Continuous Models}

In this section, starting from the description $(\ref{DIScrete})$
with a Markovian strategy and following the asymptotic
$(\ref{poiss})$ and $(\ref{diffF})$, we show that discrete
processes $(\rho_{[nt]})$ converge in distribution to solutions
of stochastic differential equations.

As in the classical case of stochastic differential equations, we
show that the evolution of a quantum system undergoing a
continuous measurement with control is either described by a
diffusive evolution or by an evolution with jump.
\begin{enumerate}
\item If $(\rho_t)$ denotes the state of a quantum system, the diffusive evolution is given by
\begin{equation}\label{eqdif}
 d\rho_t=L\big(t,u(t,\rho_t)\big)(\rho_t)dt+\Theta\big(t,u(t,\rho_t)\big)(\rho_t)dW_t
\end{equation}
where $(W_t)$ describes a one-dimensional Brownian motion. The
function $L$ is expressed as $(\ref{LINDBLAD})$ and $\Theta$ is
defined by
\begin{eqnarray}\label{deftheta}
\Theta(t,a)(\mu)=C(t,a)\mu+\mu
C^\star(t,a)-Tr\left[\mu\bigg{(}C(t,a)+C^\star(t,a)\bigg{)}\right]\mu
\end{eqnarray}
for all $t>0$, for all $a$ in $\mathbb{R}$ and all operator $\mu$
in $\mathbb{M}_2(\mathbb{C})$.
\item The evolution with jump is given by
\begin{eqnarray}\label{jjjjjj}
 d\rho_t&=&L\big(t,u(t,\rho_t)\big)(\rho_t)dt\nonumber\\&&+\left[\frac{\mathcal{J}
 \big(t,u(t,\rho_t)\big)(\rho_t)}{Tr\Big[\mathcal{J}\big(t,u(t,\rho_t)\big)(\rho_t)\Big]}
 -\rho_t\right]\big{(}d\tilde{N}_t-Tr\Big[\mathcal{J}\big(t,u(t,\rho_t)\big)(\rho_t)\Big]dt\big{)}
\end{eqnarray}
where $\tilde{N}_t$ is a counting process with stochastic
intensity
$\int_0^tTr\Big[\mathcal{J}\big(s,u(s,\rho_s)\big)(\rho_s)\Big]ds$.
The functions $L$ and $\mathcal{J}$ are as $(\ref{LINDBLAD})$.
\end{enumerate}
In a natural way, we call such equations \textit{"Controlled
Stochastic Schr\"odinger Equations"} and their solutions
\textit{"Controlled Quantum Trajectories"}.

For the moment, we do not speak about the regularity of the
functions $L$, $\Theta$ and $\mathcal{J}$. This will be discussed
when we deal with the question of existence and uniqueness of a
solution for such equations.

 This question of existence and uniqueness is relatively important because this problem is not
 really treated in details in the literature. Moreover in the two
 cases, it is difficult to show that a solution of these equations
 is valued in the set of states (actually it is a essential
 property to solve the problem of existence and uniqueness). For
 physical considerations, this property is crucial otherwise these
 equations have no sense; we are going to see that this point can
 be deduced from the convergence result. Let us start by studying the diffusive case.

\subsection{Diffusive Equation with Control}

In this section, we justify the diffusive model
$$d\rho_t=L\big(t,u(t,\rho_t)\big)(\rho_t)dt+\Theta\big(t,u(t,\rho_t)\big)(\rho_t)dW_t$$
 of controlled stochastic Schr\"odinger equations by proving that the solution
 of equation $(\ref{eqdif})$ is obtained from the limit of particular
 quantum trajectories $(\rho_{[nt]})$.
 In the same time, we show that the equation $(\ref{eqdif})$ admits a
 unique solution with values in the set of states.

Start by investigating the problem of existence and uniqueness of
a solution for $(\ref{eqdif})$. For the moment, let $u$ be any
measurable function which defines a Markovian strategy as it is
expressed in Definition $2$. Usual conditions concerning
existence and uniqueness of a solution for SDE of type
$(\ref{eqdif})$ is that for all $T>0$ there exists a constant
$M(T)$ and $K(T)$ such that the function $L$ and $\Theta$ satisfy
for all $t\leq T$ and $(\mu,\rho)\in\mathbb{M}_2(\mathbb{C})^2$ :
\begin{eqnarray}\label{condition}
&&\sup\big{\{}\Vert L(t,a)(\mu)-L(t,a)(\rho)\Vert,\Vert \Theta(t,a)(\mu)-\Theta(t,a)(\rho)\Vert\big{\}}\leq K(T)\Vert\mu-\rho\Vert\nonumber\\
&&\sup\big{\{}\Vert L(t,a)(\rho)\Vert,\Vert
\Theta(t,a)(\rho)\Vert\big{\}}\leq M(T)(1+\Vert\rho\Vert+\Vert
a\Vert).
\end{eqnarray}
Such conditions is called global Lipschitz conditions. However
even in the homogeneous case without control, such conditions are
not satisfied. Indeed, in the homogeneous situation without
control, for $\Theta$ we have
$$\Theta(t,a)(\mu)=\Theta(\mu)=C\mu+\mu C^\star-Tr\left[\mu\big{(}C+C^\star\big{)}\right]\mu.$$
Such function is not Lipschitz. Nevertheless it is $C^\infty$ and
then local Lipschitz. Such property is used in the classical case
to obtain the existence and the uniqueness of a solution for
stochastic Schr\"odinger equations (see \cite{Poisson} and
\cite{diffusion}). In the non-homogeneous context with control,
the local Lipschitz condition is expressed as follows. For all
integer $k>0$ and all $x\in\mathbb{R}$, define the function
$\phi^k$ by
$$\phi^k(x)=-k\mathbf{1}_{]-\infty,-k[}(x)+x\mathbf{1}_{[-k,k]}(x)+k\mathbf{1}_{]k,\infty[}(x).$$
The function $\phi^k$ is called a truncation function. Its
extension on the set of operator on $\mathbb{C}^2$ is given by
$$\tilde{\phi}^k(B)=(\phi^k(Re(B_{ij}))+i\phi^k(Im(B_{ij})))_{0\leq i,j\leq 1}.$$
Hence, the local Lipschitz condition for the functions $L$ and
$\theta$ can be expressed as follows. For all $T>0$ and for all
integer $k>0$ there exists a constant $M^k(T)$ and $K^k(T)$ such
that the function $L$ and $\Theta$ satisfy for all $t\leq T$ and
$(\mu,\rho)\in\mathbb{M}_2(\mathbb{C})^2$ :
\begin{eqnarray}\label{localcondition}
&&\Vert L(t,a)(\tilde{\phi}^k(\mu))-L(t,a)(\tilde{\phi}^k(\rho))\Vert\leq K^k(T)\Vert\mu-\rho\Vert\nonumber\\&&\Vert \Theta(t,a)(\tilde{\phi}^k(\mu))-\Theta(t,a)(\tilde{\phi}^k(\rho))\Vert\leq K^k(T)\Vert\mu-\rho\Vert\nonumber\\
&&\sup\big{\{}\Vert L(t,a)(\tilde{\phi}^k(\rho))\Vert,\Vert
\Theta(t,a)(\tilde{\phi}^k(\rho))\Vert\big{\}}\leq
M^k(T)(1+\Vert\rho\Vert+\Vert a\Vert).
\end{eqnarray}
As a consequence, we have the following existence and uniqueness
theorem.
\begin{thm}\label{exist} Let $u$ be any measurable function. Let $k>0$ be an integer. Let $(\Omega,\mathcal{F},\mathcal{F}_t,P)$ be a probability space which supports a standard Brownian motion $(W_t)$. Assume that $L$ and $\Theta$ satisfy the conditions $(\ref{localcondition})$. Let $\rho_0$ be any $2\times 2$ matrix. The stochastic differential equation
 \begin{equation}\label{trunc}
\rho^{\mathbf{u},k}_t=\rho_0+\int_0^tL\big(s,u(s,\tilde{\phi}^k(\rho^{\mathbf{u},k}_s)\big)
(\tilde{\phi}^k(\rho^{\mathbf{u},k}_s))ds
+\int_0^t\Theta\big(s,u(s,\tilde{\phi}^k(\rho^{\mathbf{u},k}_s))\big)(\tilde{\phi}^k(\rho^{\mathbf{u},k}_s))dW_s,
\end{equation}
admits a unique solution $(\rho^{\mathbf{u},k}_t)$. Furthermore
the application $t\mapsto\rho^{\mathbf{u},k}_t$ is almost surely
continuous.
\end{thm}

This theorem is just a consequence of the local Lipschitz
condition $(\ref{localcondition})$ (cf \cite{MR2020294}). The
process $(\rho^{\mathbf{u},k}_t)$ is called a truncated solution.
The link between such solution and a solution of the equation
$(\ref{eqdif})$ without truncature is expressed as follows.
Usually, we define the random stopping time
$$T_k=\inf\{t>0/\exists(ij),\,Re(\rho^{\mathbf{u},k}_t(ij))=k\,\,or\,\,Im(\rho^{\mathbf{u},k}_t(ij))=k\}$$
For any $k>1$, we have $T_k>0$ almost surely for $\rho_0$ is a
state and the almost surely continuity of
$(\rho^{\mathbf{u},k}_t)$ (the coefficients of $\rho_0$ satisfy
namely $\vert \rho_0(ij)\vert\leq1$). Furthermore on $[0,T_k[$ we
have
$$\tilde{\phi}^k(\rho^{\mathbf{u},k}_t)=\rho^{\mathbf{u},k}_t.$$
Therefore the process $(\rho^{\mathbf{u},k}_t)$ satisfies on
$[0,T_k[$
\begin{equation}
\rho^{\mathbf{u},k}_t=\rho_0+\int_0^tL\big(s,u(s,\rho^{\mathbf{u},k}_s)\big)
(\rho^{\mathbf{u},k}_s)ds
+\int_0^t\Theta\big(s,u(s,\rho^{\mathbf{u},k}_s)\big)(\rho^{\mathbf{u},k}_s)dW_s,
\end{equation}
Hence the process $(\rho^{\mathbf{u},k}_t)$ solution of
$(\ref{trunc})$ is the unique solution of the equation
$(\ref{eqdif})$ on $[0,T_k[$.

 In our situation, we will prove that $T_k=\infty$ for all $k>1$ by proving that the process $(\rho^{\mathbf{u},k}_t)$ is valued in the set of states. Indeed if the process $(\rho^{\mathbf{u},k}_t)$ takes value in the set of states, we have for all $t\geq0$
$$\tilde{\phi}^k(\rho^{\mathbf{u},k}_t)=\rho^{\mathbf{u},k}_t,$$
then $T_k=\infty$. As a consequence the process
$(\rho^{\mathbf{u},k}_t)$ satisfies for all $t>0$ the equation
$(\ref{eqdif})$. The truncature method becomes actually not
necessary, it just allow to exhibit a solution. As a consequence,
we have to prove that the solution obtained with a truncature
method takes value in the set of states. This property follows
from the convergence theorem.

 Indeed, let assume that there is a discrete quantum trajectory $(\rho^\mathbf{u}_{[nt]})$ which converges in distribution to $(\rho^{\mathbf{u},k}_t)$ (for some $k>1$). Such convergence is denoted by
$$\rho^\mathbf{u}_{[nt]}\Longrightarrow \rho^{\mathbf{u},k}_t.$$
Therefore for all measurable functions $\mathcal{V}$ defined on
$\mathbb{M}_2(\mathbb{C})$, we have
$$\mathcal{V}(\rho^\mathbf{u}_{[nt]})\Longrightarrow \mathcal{V}(\rho^{\mathbf{u},k}_t)$$
We apply it for the functions $\mathcal{V}(\rho)=Tr[\rho]$, for
$\mathcal{V}(\rho)=\rho^\star-\rho$ and
$\mathcal{V}_z(\rho)=\langle z,\rho z\rangle$ for all
$z\in\mathbb{C}^2$. By definition if $\rho$ is a state we have
from trace property $Tr[\rho]=1$, from self-adjointness
$\rho^\star-\rho=0$ and from positivity $\langle z,\rho
z\rangle\geq0$ for all $z\in\mathbb{C}^2$. As discrete quantum
trajectories take values in the set of states, these properties
are then conserved at the limit. The limit process
$\rho^{\mathbf{u},k}_t$ takes then also values in the set of
states. In \cite{diffusion}, the problem of existence and
uniqueness is proved independently of the convergence result. In
this case, using convergence result is more practical because the
equation is more complicated. Let us prove now the convergence
result.

Back to the description $(\ref{DIScrete})$ of discrete quantum
trajectories, with asymptotic $(\ref{diffF})$ in the case of a
non-diagonal observable $A$ and with a Markovian strategy, we have
\begin{eqnarray}\label{difff}
\rho_{[nt]}^\mathbf{u}&=&\rho_0+\sum_{k=1}^{[nt]-1}\frac{1}{n}\Big{[}L\big(k/n,u(k/n,\rho_k^\mathbf{u})\big)(\rho^\mathbf{u}_k)+\circ\left(1\right)\Big{]}\\&&+\sum_{k=1}^{[nt]-1}
\Big{[}\Theta\big(k/n,u(k/n,\rho_k^\mathbf{u})\big)(\rho_k^\mathbf{u})+\circ
(1)\Big{]}\frac{1}{\sqrt{n}}X_{k+1}(n).\nonumber
\end{eqnarray}
From this description, we can define the following processes and
functions:
\begin{eqnarray}\label{prO}
W_n(t)&=&\frac{1}{\sqrt{n}}\sum_{k=1}^{[nt]}X_{k}(n)\nonumber\\ V_n(t)&=&\frac{[nt]}{n}\nonumber\\
\rho_n^\mathbf{u}(t)&=&\rho^\mathbf{u}_{[nt]}(n)\nonumber\\
u_n(t,W)&=&u([nt]/n,W)\,\,\nonumber\\
\Theta_n(t,s)&=&\Theta([nt]/n,s)\,\,\nonumber\\
L_n(t,s)&=&L([nt]/n,s)\,\,
\end{eqnarray}
for all $t>0$, for all $s\in\mathbb{R}$ and for all
$W\in\mathbb{M}_2(\mathbb{C})$.

 By observing that these processes and
these functions are piecewise constant, we can describe the
discrete quantum trajectory $(\rho_n^\mathbf{u}(t))$ as a solution
of the following stochastic differential equation
\begin{eqnarray}\label{discretesto}
\rho^\mathbf{u}_n(t)&=&\rho_0+\int_{0}^{t}\left[L_n\Big(s\mbox{\tiny{$-$}},u_n(s\mbox{\tiny{$-$}},\rho^\mathbf{u}_n(s\mbox{\tiny{$-$}}))\Big)(\rho^\mathbf{u}_n(s\mbox{\tiny{$-$}}))+\circ(1)\right]dV_n(s)
\nonumber\\&&+\int_{0}^{t}\left[\Theta_n\Big(s\mbox{\tiny{$-$}},u_n(s\mbox{\tiny{$-$}},\rho^\mathbf{u}_n(s\mbox{\tiny{$-$}}))\Big)(\rho^\mathbf{u}_n(s\mbox{\tiny{$-$}}))+\circ(1)\right]dW_n(s)\nonumber\\
&=&\rho_0+\int_{0}^{t}\left[L_n\Big(s\mbox{\tiny{$-$}},u_n(s\mbox{\tiny{$-$}},\tilde{\phi}^k(\rho^\mathbf{u}_n(s\mbox{\tiny{$-$}}))\Big)(\tilde{\phi}^k(\rho^\mathbf{u}_n(s\mbox{\tiny{$-$}})))+\circ(1)\right]dV_n(s)
\nonumber\\&&+\int_{0}^{t}\left[\Theta_n\left(s\mbox{\tiny{$-$}},u_n(s\mbox{\tiny{$-$}},\tilde{\phi}^k(\rho^\mathbf{u}_n(s\mbox{\tiny{$-$}}))\right)(\tilde{\phi}^k(\rho^\mathbf{u}_n(s\mbox{\tiny{$-$}})))+\circ(1)\right]dW_n(s)
\end{eqnarray}
for all $k>1$. It appears essentially as a discrete version of
equation $(\ref{eqdif})$. This procedure was also used in
\cite{diffusion}, but the equation without control is more simple
and the convergence result needs less assumptions and arguments.

Let us present the arguments in the control framework. In order to
prove the convergence of this process to the solution of the
equation $(\ref{trunc})$ given by Theorem $\ref{exist}$, we use a
results of Kurtz and Protter (\cite{MR1112406}, \cite{MR1119837})
concerning weak convergence of stochastic integrals. Let us fix
some notations.

 For all $T>0$ we define $\mathcal{D}[0,T]$ the space of c\`adl\`ag process of
  $\mathbb{M}_2(\mathbb{C})$ endowed with the Skorohod topology.

 Let $T_1[0,\infty)$ denote the set of non decreasing mapping $\lambda$ from $[0,\infty)$ to $[0,\infty)$ with $\lambda(0)=0$ such that $\lambda(t+h)-\lambda(t)\leq h$ for all $t,h\geq 0$. For any function $G$ defined from $\mathbb{R}^+\times\mathbb{M}_2(\mathbb{C})$ to $\mathbb{M}_2(\mathbb{C})$, we define
$$\begin{array}{cccc}
 \tilde{G}:&\mathcal{D}[0,\infty)\times T_1[0,\infty)&\longrightarrow&\mathcal{D}[0,\infty)\\
&(X,\lambda)&\longmapsto&G(X)\circ\lambda,
\end{array}$$
such that for all $t\geq0$ we have
$G(X)\circ\lambda(t)=G(\lambda(t),X_{\lambda(t)})$. We consider
the same definition for all other functions. We introduce the two
following conditions concerning a function $\tilde{G}$ and a
sequence $\tilde{G}_n$ as above.
\begin{eqnarray}
 (C1)&& \textrm{For each compact subset}\,\,\mathcal{K}\in\mathcal{D}[0,\infty)\times T_1[0,\infty)\,\,\textrm{and}\,\,t>0,\nonumber\\
&&\sup_{(X,\lambda)}\sup_{s\leq t}\Vert\tilde{G}_n(X,\lambda)(s)-\tilde{G}(X,\lambda)(s)\Vert\rightarrow0.\nonumber\\
(C2)&&\textrm{For}\,\,(X_n,\lambda_n)_n\in\mathcal{D}[0,\infty)\times T_1[0,\infty)/\sup_{s\leq T}\Vert X_n(s)-X(s)\Vert\rightarrow0\nonumber\\
&&\textrm{and}\,\,\sup_{s\leq t}\vert\lambda_n(s)-\lambda(s)\vert\rightarrow0\,\,\textrm{for each}\,\,t>0\,\,\textrm{implies}\nonumber\\
&&\sup_{s\leq
t}\Vert\tilde{G}(X_n,\lambda_n)(s)-\tilde{G}(X,\lambda)(s)\Vert\rightarrow0.
\end{eqnarray}

 Furthermore, recall that the square-bracket $[X,X]$ is defined for a semi-martingale by the
formula:
 $$[X,X]_t=X^2_t-2\int_0^tX_{s\mbox{\tiny{$-$}}}dX_s\,.$$
 We shall denote by $T_t(V)$ the total variation of a finite variation processes $V$ on the interval $[0,t]$. The Theorem of Kurtz and Protter \cite{MR1119837} that we use is the following.

\begin{thm}\label{KP} Let $(H_n,H)$ and $(K_n,K)$ be two couple of functions which satisfy the conditions $(C1)$ and $(C2)$. Let $(\mathcal{F}_t^n)$ be a filtration and let $X_n(t)$ be a $\mathcal{F}_t^n$-adapted process which satisfies
\begin{eqnarray}X_n(t)&=&X(0)+\int_{0}^{t}H_n(s,X_n(s\mbox{\tiny{$-$}}))dV_n(s)
+\int_{0}^{t}K_n(s,X_n(s\mbox{\tiny{$-$}}))dW_n(s).
\end{eqnarray}
Let $(\Omega,\mathcal{F},\mathcal{F}_t,P)$ be a probability
space. Let $X_t$ be the unique solution of
\begin{eqnarray}X_n(t)&=&X(0)+\int_{0}^{t}H(s,X_s)ds
+\int_{0}^{t}K(s,X_s)dW_s
\end{eqnarray}
where $(W_t)$ is a standard Brownian motion on
$(\Omega,\mathcal{F},\mathcal{F}_t,P)$.

Suppose that $(W_n,V_n)$ converges in distribution in the
Skorohod topology to $(W,V)$ where $V_t=t$ for all $t\geq0$ and
suppose
\begin{eqnarray}\label{good}
&&\sup_{n}\left\{\mathbf{E}^n\Big[[W_n,W_n]_t\right]\Big\}<\infty,\\
&&\,\,\,\,\,\,\sup_{n}\Big\{\mathbf{E}^n\left[T_t(V_n)\right]\Big\}<\infty.
\end{eqnarray}

 Hence the process $(X_n(t))$ converges in distribution in $\mathcal{D}[0,T]$ for all $T>0$ to the process $(X_t)$.
\end{thm}

We wish then to apply this theorem to obtain the convergence
result for discrete quantum trajectories $(\rho_{[nt]})$
described by $(\ref{discretesto})$. Concerning the convergence of
the processes $(W_n)$ and $V_n$ we use the following theorem
which is a generalization of Donsker Theorem (see
\cite{1MR837655}).

\begin{thm}\label{conv}Let $(M_n)$ be a sequence of martingales. Suppose that
\begin{eqnarray*}\mathop{\lim}_{n\rightarrow\infty}\mathbf{E}\big{[}|[M_n,M_n]_t-t|\big{]}=0.
\end{eqnarray*}
Then $M_n$ converges in distribution to a standard Brownian
motion.
\end{thm}

In our context we have the following proposition.

\begin{pr}\label{donsker} Let $(\mathcal{F}^n_t)$ be the filtration
\begin{equation}\label{filtre}
\mathcal{F}_{t}^{n}=\sigma(X_i,i\leq[nt]).
\end{equation}
The process $(W_n(t))$ defined by $(\ref{prO})$ is a
$\mathcal{F}^n_t$-martingale. We have
$$W_n(t)\Longrightarrow W_t$$
where $(W_t)$ is a standard Brownian motion.

Moreover we have
 $$\mathop{\sup}_n\mathbf{E}\Big[[W_n,W_n]_t\Big]<\infty\,.$$

Finally, we have the convergence in distribution for the process
$(W_n,V_n,)$ to $(W,V)$ when $n$ goes to infinity.
\end{pr}
\begin{pf}
Thanks to the definition of the random variable $X_i$, we have
$\mathbf{E}[X_{i+1}/\mathcal{F}^n_{i}]=0$ which implies
$\mathbf{E}\left[\sum_{i=[ns]+1}^{[nt]}X_i/\mathcal{F}^n_s\right]=0$
for $t>s$. Thus if $t>s$ we have the martingale property:
$$\mathbf{E}[W_n(t)/\mathcal{F}^n_s]=W_n(s)+\mathbf{E}\left[\frac{1}{\sqrt{n}}
\sum_{i=[ns]+1}^{[nt]}X_i/\mathcal{F}^n_s\right]=W_n(s).$$ By
definition of $[Y,Y]$ for a stochastic process we have
$$[W_n,W_n]_t=
W_n(t)^2-2\int_0^t
W_n(s\mbox{\tiny{$-$}})dW_n(s)=\frac{1}{n}\sum_{i=1}^{[nt]}X_i^2.$$
Thus we have
\begin{eqnarray*}\mathbf{E}\big{[}[W_n,W_n]_t\big{]}&=&\frac{1}{n}\sum_{i=1}^{[nt]}
\mathbf{E}[X_i^2]=\frac{1}{n}\sum_{i=1}^{[nt]}\mathbf{E}[\mathbf{E}[X_i^2/\sigma\{X_l,l<i\}]]\\&=&\frac{1}{n}\sum_{i=1}^{[nt]}1=\frac{[nt]}{n}.
\end{eqnarray*}
Hence we have
$$\sup_n\mathbf{E}\big{[}[W_n,W_n]_t\big{]}\leq t<\infty\,.$$

 Let us
prove the convergence of $(W_n)$ to $(W)$. According to Theorem
$\ref{conv}$, we must prove that
 $$\mathop{\lim}_{n\rightarrow\infty}\mathbf{E}\big{[}|[M_n,M_n]_t-t|\big{]}=0.$$
 Actually we prove a $L_2$ convergence:
$$\mathop{\lim}_{n\rightarrow\infty}\mathbf{E}\big{[}|[M_n,M_n]_t-t|^2\big{]}=0,$$
 which implies the $L_1$ convergence. In order to show this convergence, we use the following property
   $$\mathbf{E}\left[X_i^2\right]=\mathbf{E}\left[\mathbf{E}[X_i^2/\sigma\{X_l,l<i\}]\right]=1$$ and if
    $i<j$
\begin{eqnarray*}\mathbf{E}\left[(X_i^2-1)(X^2_j-1)\right]&=&\mathbf{E}\left[(X_i^2-1)(X^2_j-1)/\sigma\{X_l,l<j\}]\right]\\
    &=&\mathbf{E}\left[(X_i^2-1)\right]\mathbf{E}\left[(X^2_j-1)\right]\\&=&0\,.
\end{eqnarray*} This gives
 \begin{eqnarray*}
 \mathbf{E}\left[\left([W_n,W_n]_t-\frac{[nt]}{n}\right)^2\right]&=&\frac{1}{n^2}\sum_{i=1}^{[nt]}\mathbf{E}\left[(X_i^2-1)^2\right]
 +\frac{1}{n^2}\sum_{i<j}\mathbf{E}
 \left[(X_i^2-1)(X^2_j-1)\right]\\&=&
 \frac{1}{n^2}\sum_{i=1}^{[nt]}\mathbf{E}\left[(X_i^2-1)^2\right].
 \end{eqnarray*}
 Thanks to the fact that $p_{00}$ and $q_{00}$ are not equal to zero (because the
observable $A$ is not diagonal!) the terms
$\mathbf{E}\left[(X_i^2-1)^2\right]$ are bounded uniformly
  in $i$ so we have:
 $$\mathop{\lim}_{n\rightarrow\infty}\mathbf{E}\left[\left([W_n,W_n]_t-\frac{[nt]}{n}\right)^2\right]=0\,.$$
 As $\frac{[nt]}{n}\longrightarrow t$ in $L_2$ we have the desired
 convergence. The convergence of $(W_n,V_n)$ is then
 straightforward.
\end{pf}\\

In order to conclude to the convergence result by using Theorem
$\ref{KP}$ of Kurtz and Protter, we have to verify conditions
$(C1)$ and $(C2)$ for the functions appearing in the equation
$(\ref{discretesto})$. We consider $\tilde{L}_n$ defined by
$$\tilde{L_n}(X)\circ(\lambda)(t)=L_n(\lambda(t),u_n(\lambda(t),X_{\lambda(t)}))(X_{\lambda(t)})+\circ(1)$$
for all $t>0$, for all $\lambda\in T_1[0,\infty)$ and all
c\`adl\`ag process $(X_t)$. Let us stress that in restriction to
the processes $(\rho_t)$ which takes values in the set of states,
the $\circ$ are uniform in $(\rho_t)$, we can then consider that
the $\circ$ are uniform for all processes. We define
$\tilde{\Theta}_n$ in the same way.

With this notation, we can express the convergence theorem.

\begin{thm}\label{convergence} Let $\mathcal{F}_t^n$ be the filtration defined by $(\ref{filtre})$. Let $\rho_0$ be any state on $\mathcal{H}_0$. Let $(\rho^\mathbf{u}_{n}(t))$ be the discrete quantum trajectory satisfying:
\begin{eqnarray}
\rho^\mathbf{u}_n(t)&=&\rho_0+\int_{0}^{t}\left[L_n(s\mbox{\tiny{$-$}},u_n(s\mbox{\tiny{$-$}},\rho^\mathbf{u}_n(s\mbox{\tiny{$-$}})))(,\rho^\mathbf{u}_n(s\mbox{\tiny{$-$}}))+\circ(1)\right]dV_n(s)
\nonumber\\&&+\int_{0}^{t}\left[\Theta_n(s\mbox{\tiny{$-$}},u_n(s\mbox{\tiny{$-$}},\rho^\mathbf{u}_n(s\mbox{\tiny{$-$}})))(\rho^\mathbf{u}_n(s\mbox{\tiny{$-$}}))+\circ(1)\right]dW_n(s).
\end{eqnarray}

Let $k>1$ be any integer. Let $(\rho^{\mathbf{u},k}_t)$ be the
unique solution of
\begin{equation}\label{eds}
\rho^{\mathbf{u},k}_t=\rho_0+\int_0^tL(s,u(s,\tilde{\phi}^k(\rho^{\mathbf{u},k}_s))(\tilde{\phi}^k(\rho^{\mathbf{u},k}_s))ds
+\int_0^t\Theta(s,u(s,\tilde{\phi}^k(\rho^{\mathbf{u},k}_s)))(\tilde{\phi}^k(\rho^{\mathbf{u},k}_s))dW_s.
\end{equation}

Assume the function $u$ is sufficiently regular such that
$\tilde{L}_n$, $\tilde{\Theta}_n$, $\tilde{L}$ and
$\tilde{\Theta}$ composed with the truncature function
$\tilde{\phi}^k$ satisfy conditions $(C1)$ and $(C2)$.

Then for all $T>0$, the process $(\rho^\mathbf{u}_{n}(t))$
converges in distribution in $\mathcal{D}[0,T]$ to the process
$(\rho_t)$.

Finally the process $(\rho^\mathbf{u}_t)$ is the unique solution
of the controlled diffusive equation
\begin{equation}
\rho^{\mathbf{u}}_t=\rho_0+\int_0^tL(s,u(s,\rho^{\mathbf{u}}_s)(\rho^{\mathbf{u}}_s)ds
+\int_0^t\Theta(s,u(s,\rho^{\mathbf{u},k}_s))(\rho^{\mathbf{u}}_s)dW_s.
\end{equation}
\end{thm}

\begin{pf}
As the condition $(C1)$ and $(C2)$ are assumed to be satisfied,
thanks to Proposition $\ref{donsker}$ and Theorem $\ref{KP}$, we
have the convergence result. The final part of the theorem comes
from the fact that the property of being a state is conserved by
passage to the limit (see the remark at the beginning of this
section).
\end{pf}\\

As regards conditions $(C1)$ and $(C2)$, the assumption for the
function $u$ is satisfied for example when $u$ is continuous. By
definition of the functions $L_n$ and $\Theta_n$ conditions
$(C1)$ and $(C2)$ are namely satisfied for the functions $L$ and
$\Theta$ satisfy the local Lipschitz conditions
$(\ref{localcondition})$ (used in Theorem $\ref{exist}$ of
existence and uniqueness).

Hence, the model of diffusive stochastic differential equation
$(\ref{eqdif})$ for continuous measurement with control is
physically justified by proving that solutions of such equations
are obtained by limit of concrete discrete procedures. In the next
section, we show a similar result by considering continuous limit
of discrete quantum trajectories of type $(\ref{poiss})$.

\subsection{Jump Equation with Control}

In this section, we investigate the convergence of a discrete
quantum trajectory which comes from repeated measurements of a
diagonal observable.

 In all this section, we fix a strategy $\mathbf{u}$ which defines a Markovian strategy.
  Furthermore, as in the diffusive case, we suppose that this strategy is continuous. Let $A$ be any diagonal observable. With the use of description $(\ref{poiss})$ and $(\ref{DIScrete})$, the discrete quantum trajectory satisfies
\begin{eqnarray}\label{PPoi}
 \rho^\mathbf{u}_{[nt]}&=&\rho_0+\sum_{k=0}^{[nt]-1}\frac{1}{n}\Big{[}L\big(
 k/n,u(k/n,\rho^\mathbf{u}_k)\big)(\rho^\mathbf{u}_k)
 -\mathcal{J}\big(k/n,u(k/n,\rho^\mathbf{u}_k)\big)(\rho^\mathbf{u}_k)\nonumber
 \\&&+Tr\big[\mathcal{J}\big(k/n,u(k/n,\rho^\mathbf{u}_k)\big)
 (\rho^\mathbf{u}_k)\big]\,\rho^\mathbf{u}_k+\circ(1)\Big{]}\nonumber\\&&
 +\sum_{k=0}^{[nt]-1}\left[\frac{\mathcal{J}\big(k/n,u(k/n,\rho^\mathbf{u}_k)\big)
(\rho^\mathbf{u}_k)}{Tr\big[\mathcal{J}(k/n,u(k/n,\rho^\mathbf{u}_k))(\rho^\mathbf{u}_k)\big]}
-\rho^\mathbf{u}_k+\circ(1)\right]\nu_{k+1}.\nonumber\\
\end{eqnarray}

Following the idea presented in article \cite{Poisson}, we aim to
show that the process $(\rho_{[nt]})$ converges
$(n\rightarrow\infty)$ to a process $(\rho_t)$
 which satisfies
\begin{eqnarray}\label{eqp1}
 \rho^\mathbf{u}_t&=&\rho_0+\int_0^t
 \bigg{[}L\big(s\mbox{\tiny{$-$}},u(s\mbox{\tiny{$-$}},\rho^\mathbf{u}_{s\mbox{\tiny{$-$}}})\big)
 (\rho^\mathbf{u}_{s\mbox{\tiny{$-$}}})+Tr\big[
 \mathcal{J}\big(s\mbox{\tiny{$-$}},u(s\mbox{\tiny{$-$}},
 \rho^\mathbf{u}_{s\mbox{\tiny{$-$}}})\big)
 (\rho^\mathbf{u}_{s\mbox{\tiny{$-$}}})\big]\,
 \rho^\mathbf{u}_{s\mbox{\tiny{$-$}}}\nonumber\\&&\hspace{2cm}-\mathcal{J}\big(
 s\mbox{\tiny{$-$}},u(s\mbox{\tiny{$-$}},\rho^\mathbf{u}_{s\mbox{\tiny{$-$}}})\big)
 (\rho^\mathbf{u}_{s\mbox{\tiny{$-$}}})\bigg{]}ds\\&&+\int_0^t\int_{\mathbb{R}}\left[
 \frac{\mathcal{J}\big(s\mbox{\tiny{$-$}},u(s\mbox{\tiny{$-$}}
 ,\rho^\mathbf{u}_{s\mbox{\tiny{$-$}}})\big)
 (\rho^\mathbf{u}_{s\mbox{\tiny{$-$}}})}{Tr\big[\mathcal{J}\big(
 s\mbox{\tiny{$-$}},u(s\mbox{\tiny{$-$}},\rho^\mathbf{u}_{s\mbox{\tiny{$-$}}})\big)
 (\rho^\mathbf{u}_{s\mbox{\tiny{$-$}}})\big]}-\rho^\mathbf{u}_{s\mbox{\tiny{$-$}}}\right]\mathbf{1}_{0<x<Tr[\mathcal{J}(s\mbox{\tiny{$-$}},u(s\mbox{\tiny{$-$}},\rho^\mathbf{u}_{s\mbox{\tiny{$-$}}}))(\rho^\mathbf{u}_{s\mbox{\tiny{$-$}}})]}N(ds,dx)\nonumber
\end{eqnarray}
where $N$ is a Poisson Point Process on $\mathbb{R}^2$. As a
consequence, if the process
 $(\rho^\mathbf{u}_t)$ exists, it gives rise to the process $(\tilde{N}_t)$ defined by
\begin{equation}
 \tilde{N}_t=\int_0^t\int_{\mathbb{R}}\mathbf{1}_{0<x<Tr[\mathcal{J}(s\mbox{\tiny{$-$}},u(s\mbox{\tiny{$-$}},\rho^\mathbf{u}_{s\mbox{\tiny{$-$}}}))(\rho_{s\mbox{\tiny{$-$}}})]}N(ds,dx)
\end{equation}
which is a counting process with stochastic intensity
$$t\rightarrow\int_0^tTr[\mathcal{J}(s\mbox{\tiny{$-$}},
u(s\mbox{\tiny{$-$}},\rho^\mathbf{u}_{s\mbox{\tiny{$-$}}}))(\rho^\mathbf{u}_{s\mbox{\tiny{$-$}}})]ds.$$

In \cite{Poisson}, it is shown that the use of the Poisson
process $N$ allows to give a mathematical sense to equation
(\ref{jjjjjj}) and allows to define properly the process
$(\tilde{N}_t)$. Actually in (\ref{jjjjjj}), the problem concerns
the definition and the existence of the process $(\tilde{N}_t)$.
Indeed, the stochastic intensity of this process depends on the
the solution of (\ref{jjjjjj}). In the way of writing the
equation (\ref{jjjjjj}), in order to define $\tilde{N}_t$, it
assumes implicitly the existence of the process solution whereas
the equation is driven by $(\tilde{N}_t)$.

 Now we consider the equation $(\ref{eqp1})$ as the
jump-model of continuous time measurement with control. It will
be justified later as limit of discrete quantum trajectories.

 For the moment, we deal with the problem existence
and uniqueness of a solution for this equation. Let us denote
\begin{eqnarray*}
 R(t,a)(\rho)&=&L(t,a)(\rho)+Tr[\mathcal{J}(t,a)(\rho)]\rho-\mathcal{J}(t,a)(\rho)\\
Q(t,a)(\rho)&=&\left(\frac{\mathcal{J}(t,a)(\rho)}{Tr[\mathcal{J}(t,a)(\rho)]}-\rho\right)\mathbf{1}_{Tr[\mathcal{J}(t,a)(\rho)]>0}
\end{eqnarray*}
for all $t\geq0$, for all $a\in\mathbb{R}$ and all state $\rho$.
It was obvious that $(\ref{eqp1})$ is equivalent to
 \begin{eqnarray*}
 \rho^\mathbf{u}_t&=&\rho_0+\int_0^tR\big(
 s\mbox{\tiny{$-$}},u(s\mbox{\tiny{$-$}},\rho^\mathbf{u}_{s\mbox{\tiny{$-$}}})\big)
 (\rho^\mathbf{u}_{s\mbox{\tiny{$-$}}})ds\\&&+\int_0^t\int_{\mathbb{R}}Q\big(
 s\mbox{\tiny{$-$}},u(s\mbox{\tiny{$-$}},\rho^\mathbf{u}_{s\mbox{\tiny{$-$}}})\big)
 (\rho^\mathbf{u}_{s\mbox{\tiny{$-$}}})\mathbf{1}_{0<x<Tr[\mathcal{J}(s\mbox{\tiny{$-$}},u(s\mbox{\tiny{$-$}},\rho^\mathbf{u}_{s\mbox{\tiny{$-$}}}))
 (\rho_{s\mbox{\tiny{$-$}}})]}N(ds,dx).\nonumber
\end{eqnarray*}

In order to prove existence and uniqueness of a solution for such
equations, a sufficient condition concerns the Lipschitz property
for functions $R$ an $\mathcal{J}$. In the same fashion of the
diffusive case, this is not the case and a truncature method is
used again. In the same way, we will have next to prove that the
truncated solution takes values in the set
of states.  Regarding functions $R$ an $\mathcal{J}$, the conditions for the Poisson case are expressed in the following remark.\\

\textbf{Remark}: As in the diffusive case, this remark concerns
the regularity of the different functions. Firstly we suppose
that $R$ and $\mathcal{J}$ satisfy the local Lipschitz condition
$(\ref{localcondition})$ defined in Section $2.1$. Secondly as
the set of states is compact, we can suppose for the stochastic
intensity that for all $T>0$ there exists a constant $K(T)$ such
that
$$Tr[\mathcal{J}(t,u(t,X_t))(X_t)]\leq K(T)$$
for all $t\geq T$ and for all c\`adl\`ag process $(X_t)$ with
values in $\mathbb{M}_2(\mathbb{C})$. This previous condition
implies the fact the stochastic intensity is bounded. Finally in
order to consider the stochastic differential equation for all
c\`adl\`ag process, we consider the function
\begin{equation}
 \tilde{Q}(t,a)(\rho)=\left(\frac{\mathcal{J}(t,a)(\rho)}
 {Re\big{(}Tr[\mathcal{J}(t,a)(\rho)]\big{)}}-\rho\right)\mathbf{1}_{Re(Tr[\mathcal{J}(t,a)(\rho)])>0}
\end{equation}
and the stochastic differential equation
\begin{eqnarray}\label{eqp2}
 &&\rho^{\mathbf{u},k}_t=\rho_0+\int_0^tR\big(
 s\mbox{\tiny{$-$}},u(s\mbox{\tiny{$-$}},\tilde{\phi}^k
 (\rho^{\mathbf{u},k}_{s\mbox{\tiny{$-$}}})\big)
 (\tilde{\phi}^k(\rho^{\mathbf{u},k}_{s\mbox{\tiny{$-$}}}))ds\\&&+
 \int_0^t\int_{\mathbb{R}}\tilde{Q}\big(
 s\mbox{\tiny{$-$}},u(s\mbox{\tiny{$-$}},
 \tilde{\phi}^k(\rho^{\mathbf{u},k}_{s\mbox{\tiny{$-$}}})\big)
 (\tilde{\phi}^k(\rho^{\mathbf{u},k}_{s\mbox{\tiny{$-$}}}))\mathbf{1}_{0<x<Re(Tr[\mathcal{J}(s\mbox{\tiny{$-$}},u(s\mbox{\tiny{$-$}},\tilde{\phi}^k(\rho^{\mathbf{u},k}_{s\mbox{\tiny{$-$}}})))(\tilde{\phi}^k(\rho^{\mathbf{u},k}_{s\mbox{\tiny{$-$}}}))]}N(ds,dx).\nonumber
\end{eqnarray}
where $\tilde{\phi}^k$ is a truncature function defined in Section
$3.1$. As in the diffusive case, if a solution of the equation
$(\ref{eqp2})$ takes value in the set of states, it is a solution
of the equation $(\ref{eqp1})$. In addition to the diffusive case,
we have to remark that if $\rho$ is a state
$$Re(Tr[\mathcal{J}(t,a)(\rho)])=Tr[\mathcal{J}(t,a)(\rho)]\geq0$$
for all $t\geq0$ and for all $a\in\mathbb{R}$.\\

Exactly in the same way as the diffusive case, if we show that a
discrete quantum trajectory converges in distribution to a
solution of the truncated equation $(\ref{eqp2})$, it involves
that this solution takes values in the set of states. Let us
first deal with the problem of existence and uniqueness of a
solution for the equation $(\ref{eqp2})$. We have the following
theorem due to Jacod and Protter in \cite{QR}.

\begin{thm} Let $(\Omega,\mathcal{F},P)$ be a probability space of a Poisson point Process $N$. The stochastic differential equation
\begin{eqnarray}
 &&\rho^{\mathbf{u},k}_t=\rho_0+\int_0^tR\big(
 s\mbox{\tiny{$-$}},u(s\mbox{\tiny{$-$}},\tilde{\phi}^k
 (\rho^{\mathbf{u},k}_{s\mbox{\tiny{$-$}}})\big)
 (\tilde{\phi}^k(\rho^{\mathbf{u},k}_{s\mbox{\tiny{$-$}}}))ds\\&&+
 \int_0^t\int_{\mathbb{R}}\tilde{Q}\big(
 s\mbox{\tiny{$-$}},u(s\mbox{\tiny{$-$}},
 \tilde{\phi}^k(\rho^{\mathbf{u},k}_{s\mbox{\tiny{$-$}}})\big)
 (\tilde{\phi}^k(\rho^{\mathbf{u},k}_{s\mbox{\tiny{$-$}}}))\mathbf{1}_{0<x<Re(Tr[\mathcal{J}(s\mbox{\tiny{$-$}},u(s\mbox{\tiny{$-$}},\tilde{\phi}^k(\rho^{\mathbf{u},k}_{s\mbox{\tiny{$-$}}})))(\tilde{\phi}^k(\rho^{\mathbf{u},k}_{s\mbox{\tiny{$-$}}}))]}N(ds,dx)\nonumber
\end{eqnarray}
 admits a unique solution $\rho^{\mathbf{u},k}_t$ defined for al$t\geq0$. Furthermore the process
$$\overline{N}_t=\int_0^t\int_{\mathbb{R}}\mathbf{1}_{0<x<Re(Tr[\mathcal{J}(s\mbox{\tiny{$-$}},u(s\mbox{\tiny{$-$}},\tilde{\phi}^k(\rho^{\mathbf{u},k}_{s\mbox{\tiny{$-$}}})))(\tilde{\phi}^k(\rho^{\mathbf{u},k}_{s\mbox{\tiny{$-$}}}))]}N(ds,dx)$$
allows to define the filtration $(\overline{\mathcal{F}}_t)$ where
$\overline{\mathcal{F}}_t=\sigma\{\overline{N}_s,s\leq t\}.$

Hence the process
$$\overline{N}_t-\int_0^t\bigg{[}Re(Tr[\mathcal{J}(s\mbox{\tiny{$-$}},u(s\mbox{\tiny{$-$}},\tilde{\phi}^k(\rho^{\mathbf{u},k}_{s\mbox{\tiny{$-$}}})))(\tilde{\phi}^k(\rho^{\mathbf{u},k}_{s\mbox{\tiny{$-$}}}))\bigg{]}^+ds$$
is a $\overline{\mathcal{F}}_t$-martingale.
\end{thm}

The term $(x)^+$ denotes $\max(0,x)$. Such theorem is treated in
details in \cite{Poisson} for quantum trajectories without
control. We give here a way to express the solution of
$(\ref{eqp2})$ in a particular case.

 Suppose that there exists a
constant $K$ such that:
\begin{equation}\label{con}
\bigg{[}Re(\mathcal{J}(t,u(t,X_t))(X_t))\bigg{]}^+<K,
\end{equation} for all $t\geq0$ and all c\`adl\`ag process $(X_t)$. With this property we can consider only the points of $N$ contained in $\mathbb{R}\times[0,K]$. The random
function $$\mathcal{N}_t:t\rightarrow N(.,[0,t]\times[0,K])$$
defines then a standard Poisson process with intensity $K$. Let
$T>0$, the Poisson Random Measure and the previous process
generate on $[0,T]$ a sequence
$\{(\tau_i,\xi_i),i\in\{1,\ldots,\mathcal{N}_t)\}\}$. Each
$\tau_i$ represents the jump time of the process
$(\mathcal{N}_t)$. Moreover the random variables $\xi_i$ are
random uniform variables on $[0,K]$. Let $k>1$ be a fixed integer,
we can write the solution of $(\ref{eqp2})$ in the following way:
 \begin{eqnarray}\label{expression}
\rho^{\mathbf{u},k}_{t}&=&\rho_0+\int_0^tR(s\mbox{\tiny{$-$}},u(s\mbox{\tiny{$-$}},\tilde{\phi}^k(\rho^{\mathbf{u},k}_{s\mbox{\tiny{$-$}}})))
(\tilde{\phi}^k(\rho^{\mathbf{u},k}_{s\mbox{\tiny{$-$}}}))ds
\nonumber\\&&+\sum_{i=1}^{\mathcal{N}_t}Q(\tau_{i}-,u(\tau_{i}-,\tilde{\phi}^k(\rho^{\mathbf{u},k}_{\tau_i-}))
(\tilde{\phi}^k(\rho^{\mathbf{u},k}_{\tau_i-})) \mathbf{1}_{0\leq
\xi_i\leq
(Re(Tr[\mathcal{J}(\tau_i-,u(\tau_i-,\tilde{\phi}^k(\rho^{\mathbf{u},k}_{\tau_i-})))
(\tilde{\phi}^k(\rho^{\mathbf{u},k}_{\tau_i-})]))^+}\nonumber\\
\overline{N}_t&=&\sum_{i=1}^{\mathcal{N}_t}\mathbf{1}_{0\leq
\xi_i\leq
(Re(Tr[\mathcal{J}(\tau_i-,u(\tau_i-,\tilde{\phi}^k(\rho^{\mathbf{u},k}_{\tau_i-})))
(\tilde{\phi}^k(\rho^{\mathbf{u},k}_{\tau_i-})]))^+}.
 \end{eqnarray}
The general case is treated in details in \cite{QR}. Let us make
more precise how the solution of $(\ref{eqp2})$ is defined from
the expression $(\ref{expression})$ in the particular case
$(\ref{con})$. By applying Cauchy-Lipschitz Theorem, we consider
the solution of the ordinary differential equation
\begin{equation}\label{ordi}\rho^{\mathbf{u},k}_{t}=\rho_0+\int_0^tR(s\mbox{\tiny{$-$}},u(s\mbox{\tiny{$-$}},\tilde{\phi}^k(\rho^{\mathbf{u},k}_{s\mbox{\tiny{$-$}}})))
(\tilde{\phi}^k(\rho^{\mathbf{u},k}_{s\mbox{\tiny{$-$}}}))ds.\end{equation}
It gives rise to the function
$$t\mapsto\bigg{[}Re(Tr[\mathcal{J}(t,u(t,\tilde{\phi}^k(\rho^{\mathbf{u},k}_{t})))
(\tilde{\phi}^k(\rho^{\mathbf{u},k}_{t}))\bigg{]}^+.$$ Let define
the first jump-time of the process $(\overline{N}_t)$. For this,
we introduce the set
$$G_t=\{(x,y)\in\mathbb{R}^2/0<x\leq t,0<y<\bigg{[}Re(Tr[\mathcal{J}(x,u(x,\tilde{\phi}^k(\rho^{\mathbf{u},k}_{x})))
(\tilde{\phi}^k(\rho^{\mathbf{u},k}_{x}))\bigg{]}^+\},$$ and the
random stopping time
$$T_1=\inf\{t/N(G_t)=1\}.$$
As a consequence on $[0,T_1[$ the solution of $(\ref{eqp2})$ is
given by the solution of the ordinary differential equation
(\ref{ordi}) and $\rho^{\mathbf{u},k}_{T_1}$ is defined by
$$\rho^{\mathbf{u},k}_{T_1}=\rho^{\mathbf{u},k}_{T_{1}-}+
Q(T_1-,u(T_1-,\rho^{\mathbf{u},k}_{T_{1}-}))(\rho^{\mathbf{u},k}_{T_{1}-}).$$
We solve the ordinary differential equation after $T_1$ with this
initial condition, and by a similar reasoning we shall define a
second jump-time. Thus, we construct a sequence of jump-time
$(T_n)$. The boundness property $(\ref{con})$ implies that the
stochastic intensity is bounded. Hence, we can show $\lim
T_n=\infty$ almost surely (see \cite{QR} or \cite{Poisson}).

The solution of $(\ref{eqp2})$ is then given by the solution of
the ordinary differential equation
$$d\rho^{\mathbf{u},k}_{t}=R(t,u(t,\tilde{\phi}^k(\rho^{\mathbf{u},k}_{t})))
(\tilde{\phi}^k(\rho^{\mathbf{u},k}_{t}))dt$$ between the jump of
the process $\overline{N}_t$. The process $\overline{N}_t$
corresponds to the number of point of the Poisson point process
$N$ included in the x axis and the curve
$$t\mapsto\bigg{[}Re(Tr[\mathcal{J}(t,u(t,\tilde{\phi}^k(\rho^{\mathbf{u},k}_{t})))
(\tilde{\phi}^k(\rho^{\mathbf{u},k}_{t}))\bigg{]}^+.$$ The general
case is more technical but can be expressed in the same way.
\bigskip

In order to summarize the procedure to show the existence and
uniqueness of solution, it is worth noticing that all the
technical precaution are justified because we do not know that
the equation preserves the property of being a state before
showing the convergence result (it is the same problem in the
diffusive case).
\smallskip

 Now we
investigate the convergence result. In a first time, the way to
proceed is the same as in \cite{Poisson}. Next, we use another
way to obtain the convergence result because we cannot applied
the Theorem of Kurtz and Protter in this case.

From expression $(\ref{PPoi})$, define
\begin{eqnarray*}
\rho^\mathbf{u}_n(t)&=&\rho^\mathbf{u}_{[nt]}\\
 N_n(t)&=&\sum_{k=1}^{[nt]}\nu_{k},\,\,\\V_n(t)&=&\frac{[nt]}{n}\\
R_n(t,a)(\rho)&=&R([nt]/n,a)(\rho),\,\,\\
Q_n(t,a)(\rho)&=&Q([nt]/n,a)(\rho)\\
u_n(t,W)&=&u([nt]/n,W)
\end{eqnarray*}
for all $t\geq0$, for all $a\in\mathbb{R}$ and all
$W\in\mathbb{M}_2(\mathbb{C})$. Hence the process
$(\rho^\mathbf{u}_n(t))$ satisfies the stochastic differential
equation
\begin{eqnarray*}
\rho^\mathbf{u}_n(t)&=&\int_{0}^t\left[R_n\big(
s\mbox{\tiny{$-$}},
u_n(s\mbox{\tiny{$-$}},\rho^\mathbf{u}_n(s\mbox{\tiny{$-$}})\big)
(\rho^\mathbf{u}_n(s\mbox{\tiny{$-$}}))+\circ(1)\right]dV_n(s)\\
&&+\int_{0}^t\left[Q_n\big(s
\mbox{\tiny{$-$}},u_n(s\mbox{\tiny{$-$}}
,\rho^\mathbf{u}_n(s\mbox{\tiny{$-$}}))(\rho^\mathbf{u}_n(s\mbox{\tiny{$-$}})\big)+\circ(1)\right]dN_n(s).
\end{eqnarray*}

In this case, we do not have directly an equivalent of the Donsker
Theorem for the process $(N_n(t))$. Because of the stochastic
intensity of $(\tilde{N}_t)$ which depends on the solution, it is
actually not possible to prove the convergence of $(N_n(t))$ to
$(\tilde{N}_t)$ independently of the convergence of
$(\rho_{[nt]})$ to $(\rho_t)$. Hence we cannot applied the result
of Kurtz and Protter used in the diffusive case.

 The convergence result is here
obtained by using a random coupling method, that is, we realize
the process $(\rho_{[nt]})$ in the probability space of the
Poisson Point Process $N$. This method allows then to compare
directly continuous and discrete quantum trajectories in the same
probability space. It is described as follows.

Remember that the random variables $(\mathbf{1}_1^{k})$ satisfy:
\begin{equation*}\left\{\begin{array}{ccl}
\mathbf{1}^{k+1}_1(0)=0&\textrm{with
probability}&p_{k+1}(n)=1-\frac{1}{n}Tr\big[
\mathcal{J}(k/n,u(k/n,\rho^\mathbf{u}_k))(\rho^\mathbf{u}_k)\big]+\circ\left(\frac{1}{n}\right)\\&&\\
\mathbf{1}^{k+1}_1(1)=1&\textrm{with
probability}&q_{k+1}(n)=\frac{1}{n}Tr\big[
\mathcal{J}(k/n,u(k/n,\rho^\mathbf{u}_k))(\rho^\mathbf{u}_k)\big]
+\circ\left(\frac{1}{n}\right)\end{array}\right.
\end{equation*}

We define the following sequence of random variable which are
defined on the set of states
 \begin{equation}\tilde{\nu}_{k+1}(\eta,\omega)=\mathbf{1}_{N(\omega,G_k(\eta))>0}
 \end{equation}
 where $G_k(\eta)=\left\{(t,u)/\frac{k}{n}\leq t<\frac{k+1}{n},0\leq
u\leq-n\ln(Tr[\mathcal{L}^{k+1}_{0}(n)(\eta)])\right\}$. Let
$\rho_0=\rho$ be any state and $T>0$, we define the process
$(\tilde{\rho}_k)$ for $k<[nT]$ by the recursive formula
\begin{eqnarray}\label{eqq}
\tilde{\rho}^\mathbf{u}_{k+1}&=&\mathcal{L}^{k+1}_{0}(\tilde{\rho}^\mathbf{u}_k)
+\mathcal{L}^{k+1}_{1}(\tilde{\rho}^\mathbf{u}_k)\nonumber\\
&&+\left[-\frac{\mathcal{L}^{k+1}_{0}(\tilde{\rho}^\mathbf{u}_k)}{Tr[\mathcal{L}^{k+1}_{0}(\tilde{\rho}^\mathbf{u}_k)]}+
\frac{\mathcal{L}^{k+1}_{1}(\tilde{\rho}^\mathbf{u}_k)}{Tr[\mathcal{L}^{k+1}_{1}(\tilde{\rho}^\mathbf{u}_k)]}\right]
\left(\tilde{\nu}_{k+1}(\tilde{\rho}^\mathbf{u}_k,.)-Tr[\mathcal{L}^{k+1}_{1}(\tilde{\rho}^\mathbf{u}_k)]\right).
\end{eqnarray}
Thanks to properties of Poisson probability measure, the random
variables $(\mathbf{1}^{k}_1)$ and $(\tilde{\nu}_k)$ have the same
distribution. It involves the following property.

\begin{pr}Let $T$ be fixed. The discrete process $(\tilde{\rho}^\mathbf{u}_k)_{k\leq[nT]}$ defined by $(\ref{eqq})$
have the same distribution of the discrete quantum trajectory
$(\rho^\mathbf{u}_k)_{k\leq[nT]}$ defined by the quantum repeated
measurement.
\end{pr}

The convergence result is then expressed as follows.

\begin{thm}\label{convpoiss}Let $T>0$. Let $(\Omega,\mathcal{F},P)$ be a
probability space of a Poisson Point process $N$. Let
$(\tilde{\rho}^\mathbf{u}_{[nt]})_{0\leq t\leq T}$ be the discrete
quantum trajectory defined by the recursive formula $(\ref{eqq})$

Hence, for all $T>0$ the process
$(\tilde{\rho}^\mathbf{u}_{[nt]})_{0\leq t\leq T}$ converges in
distribution in $\mathcal{D}[0,T]$ (for the Skorohod topology) to
the process $(\rho_t^\mathbf{u})$ solution of the stochastic
differential equation:
\begin{eqnarray}\label{sto}
 \rho^\mathbf{u}_t&=&\rho_0+\int_{0}^tR(s\mbox{\tiny{$-$}},u(s\mbox{\tiny{$-$}},\rho^\mathbf{u}_{s\mbox{\tiny{$-$}}})(\rho^\mathbf{u}_{s\mbox{\tiny{$-$}}})ds
 \\&&+\int_{0}^t\int_{\mathbb{R}}Q(s\mbox{\tiny{$-$}},u(s\mbox{\tiny{$-$}},\rho^\mathbf{u}_{s\mbox{\tiny{$-$}}})(\rho^\mathbf{u}_{s\mbox{\tiny{$-$}}})
 \mathbf{1}_{0<x<Tr[\mathcal{J}(s\mbox{\tiny{$-$}},u(s\mbox{\tiny{$-$}},\rho^\mathbf{u}_{s\mbox{\tiny{$-$}}}))(\rho^\mathbf{u}_{s\mbox{\tiny{$-$}}})])}N(.,ds,dx).\nonumber
\end{eqnarray}
\end{thm}

This theorem relies on the fact that the process
$(\tilde{\rho}^\mathbf{u}_{[nt]})$ satisfies the same asymptotic
of the discrete quantum trajectory $(\rho^\mathbf{u}_{[nt]})$; we
have namely
\begin{eqnarray}
 \tilde{\rho}_{[nt]}&=&\rho_0+\sum_{k=0}^{[nt]-1}\frac{1}{n}\left[R(k/n,u(k/n,\tilde\rho^\mathbf{u}_k))
 (\tilde\rho^\mathbf{u}_k)+\circ(1)
 \right]\nonumber\\&&+\sum_{k=0}^{[nt]-1}\left[Q(k/n,u(k/n,\tilde\rho^\mathbf{u}_k))(\tilde\rho^\mathbf{u}_k)+
 \circ(1)\right]\tilde{\nu}_{k+1}(\rho^\mathbf{u}_k,.).
\end{eqnarray}

The complete proof of this theorem is very technical. The idea is
to compare the discrete process $(\rho_{[nt]})$ with an Euler
Scheme of the solution of the jump-equation. More details for
such techniques can be found in \cite{Poisson} where the case
without control is entirely developed.

In the next section, we expose examples and applications of such
stochastic models.

\section{Examples and Applications}

This section is devoted to some applications of quantum
measurement with control. In a first time, by a discrete model, we
justify a stochastic model for the experience of Resonance
fluorescence. The setup is the one of a laser driving an atom in
presence of a photon counter. In a second time, we present general
results in Stochastic Control theory applied to quantum
trajectories.

\subsection{Laser Monitoring Atom: Resonance Fluorescence}

We here describe a discrete model of an atom monitored by a laser.
A measurement is performed by a photon counter which detects the
photon emission. The setup of repeated quantum interactions is
described as follows.

 The length time of interaction is chosen to be $h=1/n$. Let us
 describe one interaction. Here we need three basis spaces.
  The atom system is represented by $\mathcal{H}_0$ equipped with a state $\rho$.
  The laser is representing by $(\mathcal{H}^l,\mu^l)$ and the photon counter by $(\mathcal{H}^c,\beta^c)$.
  Each Hilbert space are $\mathbb{C}^2$ endowed with the orthonormal basis $(\Omega,X)$ and
  the unitary operator is denoted by $U$. The compound system after interaction is:
$$\mathcal{H}_0\otimes\mathcal{H}^l\otimes\mathcal{H}^c,$$
and the state after interaction is:
$$\alpha=U(\rho\otimes\mu^l\otimes\beta^c)U^\star.$$
The appropriate orthonormal basis
$\mathcal{H}_0\otimes\mathcal{H}^l\otimes\mathcal{H}^c$, in this
case, is $\Omega\otimes\Omega\otimes\Omega,\,
X\otimes\Omega\otimes\Omega,\,\Omega\otimes
X\otimes\Omega,\,X\otimes
X\otimes\Omega,\Omega\otimes\Omega\otimes
X,\,X\otimes\Omega\otimes X,\,\Omega\otimes X\otimes X,\,X\otimes
X\otimes X$. As in the presentation of the discrete two level
atom in contact with a spin chain, the unitary operator is here
considered as a $4\times4$ matrix $U=(L_{i,j}(n))_{0\leq i,j\leq
3}$ where each $L_{ij}(n)$ are operator on $\mathcal{H}_0$.

If the different state of the laser and the counter as of the form
$$\mu^l=\left(\begin{array}{cc} a&b\\c&d\end{array}\right),\,\,\beta^c=\vert \Omega\rangle\langle \Omega\vert,$$
for the state $\alpha=(\alpha_{uv})_{0\leq u,v\leq 3}$ after
interaction, we have
\begin{equation}
 \alpha_{uv}=\Big(aL_{u0}(n)\rho+bL_{u1}(n)\rho\Big)L_{v0}^\star(n)+\Big(cL_{u0}(n)\rho+dL_{u1}(n)\rho\Big)L_{v1}^\star(n).
\end{equation}

The measurement is performed on the counter photon side. Let $A$
denotes any observable of $\mathcal{H}^c$ then $I\otimes I\otimes
A$ denotes the corresponding observable on
$\mathcal{H}_0\otimes\mathcal{H}^l\otimes\mathcal{H}^c.$ We
perform a measurement and by partial trace operation with respect
to $\mathcal{H}^l\otimes\mathcal{H}^c$ we obtain a new state on
$\mathcal{H}_0$.

The control is rendered by the modification at each interaction of
the intensity of the laser. This modification is here taken into
account by the reference state of the laser. The reference state
at the $k$-th interaction is denoted by $\mu_k^l$. In the
continuous case of Resonance fluorescence, the state of a laser is
usually described by a coherent vector on a Fock space (see
\cite{MR2042615}). From works of Attal and Pautrat in
approximation of Fock space (\cite{MR1971605},\cite{MR2160332}),
in our context the discrete state of the laser can be described by
\begin{equation}
 \mu_k^l=\left(\begin{array}{cc}a(k/n)&b(k/n)
 \\c(k/n)&d(k/n)\end{array}\right)=\frac{1}{1+\vert h(k/n)\vert^2}\left(\begin{array}{cc}1&h(k/n)
 \\\overline{h}(k/n)&\vert h(k/n)\vert^2\end{array}\right).
\end{equation}
The function $h$ represents the evolution of the intensity of the
laser and depends naturally on $n$.

If $\rho_k$ denotes the state on $\mathcal{H}_0$ after $k$ first
measurement,
 the state
 $$\alpha^{k+1}(n)=\left(\alpha_{uv}^{k+1}(n)\right)_{0\leq u,v\leq 3}=U_{k+1}(n)(\rho_k\otimes\mu_k^l\otimes\beta^c)U^\star_{k+1}(n)$$ after interaction
 satisfies
\begin{eqnarray*}\alpha^{k+1}_{uv}(n)&=&\Big(a(k/n)L_{u0}(n)\rho_k+b(k/n)L_{u1}(n)\rho_k\Big)L_{v0}^\star(n)\\&&+\Big(c(k/n)L_{u0}(n)\rho_k+d(k/n)L_{u1}(n)\rho_k\Big)L_{v1}^\star(n)
\end{eqnarray*}

\textbf{Remark:} Let us stress that is not directly the framework
of Section 1. Here, the control is namely not rendered by the
modification of the unitary evolution. Moreover the interacting
system is described by
$(\mathcal{H}^l\otimes\mathcal{H}^c,\mu^l_k\otimes\beta)$ and
$\mu^l_k\otimes\beta$ is not of the form $\vert X_0\rangle\langle
X_0\vert$ as in Section $1$. In order to translate this setting in
the case of discrete models of Section $1$, one can use the G.N.S
Representation theory of a finite dimensional Hilbert space
(\cite{MR1468230},\cite{MR1468229}). This theory allows to
consider the state $\mu^l_k\otimes\beta$ as a state of the form
$\vert X_0\rangle\langle X_0\vert$ in a biggest Hilbert space.
The G.N.S representation modifies then the expression of operator
$U_k$, and the control expressed in $\mu^l_k\otimes\beta$ is
again expressed in the new expression of $U_k$ (see
\cite{MR2323437} for more details). In our case, we do not use
such theory because it is more explicit to make directly
computations to reach the discrete equation in asymptotic form.
\bigskip

Let us present the result. The principle of measurement is the
same as in Section $1$. The counting case is also given by a
diagonal observable of $\mathcal{H}_c$. We shall focus on this
case which renders the emission of photon (\cite{MR2042615}). The
asymptotic for the unitary operator follows the asymptotic of
Attal-Pautrat in \cite{AP}. Let $\delta_{ij}=1$ if $i=j$ we
denote:
 $$\epsilon_{ij}=\frac{1}{2}(\delta_{0i}+\delta_{0j})$$
The coefficients must follow the convergence condition:
$$\lim_{n\rightarrow\infty}n^{\epsilon_{ij}}(L_{ij}(n)-\delta_{ij}I)=L_{ij}$$
where $L_{ij}$ are operator on $\mathcal{H}_0$.

Let $P_0=\vert \Omega\rangle\langle \Omega\vert$ and $P_1=\vert
X\rangle\langle X\vert$ be eigenprojectors of a diagonal
observable $A$. If $\rho_k$ denotes the random state after $k$
measurements we denote:
\begin{eqnarray}
 \mathcal{L}_0^{k+1}(\rho_k)&=&\mathbf{E}_0[I\otimes I\otimes P_0(U_{k+1}(n)(\rho_k\otimes\mu_k^l\otimes\beta)U^\star_{k+1}(n))I\otimes I\otimes P_0]\nonumber\\
&=&\alpha_{00}^{k+1}(n)+\alpha_{11}^{k+1}(n)\nonumber\\
\mathcal{L}_1^{k+1}(\rho_k)&=&\mathbf{E}_0[I\otimes I\otimes P_1(U_{k+1}(n)(\rho_k\otimes\mu_k^l\otimes\beta)U^\star_{k+1}(n))I\otimes I\otimes P_1]\nonumber\\
&=&\alpha_{22}^{k+1}(n)+\alpha_{33}^{k+1}(n)
\end{eqnarray}
This is namely the two non normalized state, the operator
$\mathcal{L}_0^{k+1}(\rho_k)$ appears with probability
$p_{k+1}=Tr[\mathcal{L}_0^{k+1}(\rho_k)]$ and
$\mathcal{L}_1^{k+1}(\rho_k)$ with probability
$q_{k+1}=Tr[\mathcal{L}_0^{k+1}(\rho_k)]$.

From works of Attal-Pautrat in approximation and asymptotic in
Fock space, we put
$$h(k/n)=\frac{1}{\sqrt{n}}f(k/n)+\circ\left(\frac{1}{n}\right),$$
where $f$ is a function from $\mathbb{R}$ to $\mathbb{C}$. In the
same way of Section $2$, we assume that the intensity of the laser
$f$ is continuous.

With the same arguments of Section $1$, the evolution of the
discrete quantum trajectory is described by
\begin{equation}\label{discrette}
 \rho_k=\frac{\mathcal{L}_0^{k+1}(\rho_k)}{p_{k+1}}+\left[-\frac{\mathcal{L}_0^{k+1}(\rho_k)}{p_{k+1}}+\frac{\mathcal{L}_1^{k+1}(\rho_k)}{q_{k+1}}\right]\mathbf{1}^{k+1}_1
\end{equation}

For a further use, convergence results will be established in the
case $L_{01}=-L_{10}^\star$, and $L_{11}=L_{21}=L_{31}=L_{30}=0$.
Conditions about asymptotic of $U$ and the fact that it is a
unitary-operator imply that
\begin{equation}\label{condunitary}
L_{00}=-(iH+\frac{1}{2}\sum_{i=1}^2 L_{i0}^\star L_{i0})
\end{equation}
In the same way of Section $2.2$, convergence result in this
situation is expressed as follows.

\begin{pr}
Let $(\Omega,\mathcal{F},\mathcal{F}_t,P)$ be a probability space
of a Poisson point process $N$ on $\mathbb{R}^2$.

 The discrete quantum trajectory $(\rho_{[nt]})_{0\leq t\leq T}$ defined by the equation $(\ref{discrette})$
  weakly converges in $\mathcal{D}([0,T])$ for all $T$ to the solution of the following stochastic differential equation:
\begin{eqnarray}\label{elas}
\rho_t&=&\rho_0+\int_{0}^t
\bigg{[}-i[H,\rho_{s\mbox{\tiny{$-$}}}]+\frac{1}{2}\{\sum_{i=1}^2L_{i0}^\star
L_{i0},\rho_{s\mbox{\tiny{$-$}}}\}+L_{10}\rho_{s\mbox{\tiny{$-$}}}L_{10}^\star\nonumber\\&&\hspace{1,5cm}+
[\overline{f}(s\mbox{\tiny{$-$}})L_{10}\rho_{s\mbox{\tiny{$-$}}}-f(s\mbox{\tiny{$-$}})L_{10}^\star,\rho_{s\mbox{\tiny{$-$}}}]-Tr[L_{20}\rho_{s\mbox{\tiny{$-$}}}L_{20}^\star]\rho_{s\mbox{\tiny{$-$}}}
\bigg{]}ds\nonumber\\&&+\int_{0}^t\int_{\mathbb{R}}\left[-\rho_{s\mbox{\tiny{$-$}}}+\frac{L_{20}\rho_{s\mbox{\tiny{$-$}}}L_{20}^\star}{Tr[L_{20}
\rho_{s\mbox{\tiny{$-$}}}L_{20}^\star]}\right]\mathbf{1}_{0<x<Tr[L_{20}\rho_{s\mbox{\tiny{$-$}}}L_{20}^\star]}N(dx,ds).
\end{eqnarray}
\end{pr}

\begin{pf}
 For example we have the following asymptotic for $\mathcal{L}_0^{k+1}(\rho_k)$:
\begin{eqnarray*}&&\mathcal{L}_0(\rho_k)=\\&&\rho_k+\frac{1}{n}\left[L_{00}\rho+\rho L_{00}^\star+L_{10}\rho L_{10}^\star+f(\frac{k}{n})[L_{01}\rho+\rho L_{10}^\star]+\overline{f}(\frac{k}{n})[L_{10}\rho+
\rho L_{01}^\star]\right] +\circ\left(\frac{1}{n}\right)\nonumber
\end{eqnarray*}

This above asymptotic, the condition about the operator $L_{ij}$
and the theorem $(\ref{convpoiss})$ prove the proposition.
\end{pf}

The stochastic differential equation $(\ref{elas})$ is then the
continuous time stochastic model of Resonance fluorescence.
 In this model, the control is deterministic. Before to give an application of stochastic control
 , let us briefly expose a use of the laser monitoring atom model.\\

Consider the special case, where the Hamiltonian $H=0$. Let put
$$C=\left(\begin{array}{cc}0&1\\0&0\end{array}\right),\,\,L_{10}=k_lC,\,\,L_{20}=k_cC,$$
with $\vert k_l\vert^2+\vert k_c\vert^2=1$. The constant $k_f$
and $k_c$ are called decay rates (\cite{MR2042615}).

Without control, the stochastic model of a two level atom in
presence of a photon counter (\cite{Poisson}) is given by:
\begin{eqnarray}\label{withoutcontroll}
\mu_t&=&\mu_0+\int_{0}^t
\bigg{[}+\frac{1}{2}\{C,\mu_{s\mbox{\tiny{$-$}}}\}+C\mu_{s\mbox{\tiny{$-$}}}C^\star-Tr[C\mu_{s\mbox{\tiny{$-$}}}C^\star]\mu_{s\mbox{\tiny{$-$}}}\bigg{]}ds\nonumber\\&&+\int_{0}^t\int_{\mathbb{R}}\left[-\mu_{s\mbox{\tiny{$-$}}}+\frac{C\mu_{s\mbox{\tiny{$-$}}}C^\star}{Tr[C\mu_{s\mbox{\tiny{$-$}}}C^\star]}\right]\mathbf{1}_{0<x<Tr[C\mu_{s\mbox{\tiny{$-$}}}C^\star]}N(dx,ds).
\end{eqnarray}

Let denote
$\tilde{N}_t=\int_{0}^t\int_{\mathbb{R}}\mathbf{1}_{0<x<Tr[C\mu_{s\mbox{\tiny{$-$}}}C^\star]}N(dx,ds)$
and $T=\inf\{t>0;\,\tilde{N}_t>0\}$. In \cite{commun} it was
proved that:
\begin{equation}\label{stab}
 \mu_t=\left(\begin{array}{cc}1&0\\0&0\end{array}\right)=\vert\Omega\rangle\langle\Omega\vert.
\end{equation}
for all $t>T$. Physically, it means that at most one photon
appears on the photon counter. The mathematical reason is that if
we write the equation $(\ref{withoutcontroll})$ in the following
way:
$$\mu_t=\int_0^t\Psi(\mu_{s\mbox{\tiny{$-$}}})ds+\int_0^t\Phi(\mu_{s\mbox{\tiny{$-$}}})d\tilde{N}_s,$$
we have for $\mu=\vert\Omega\rangle\langle\Omega\vert$
$$\Phi(\mu)=\Psi(\mu)=0.$$ The state
$\vert\Omega\rangle\langle\Omega\vert$ is an equilibrium state.

In the presence of laser, the control $f$ gives rise to the term
 $[\overline{f}L_{10}-fL_{10}^\star,.]=[k_l\overline{f}C-\overline{k_l}fC^\star,.]$. Hence if $\mu=\vert\Omega\rangle\langle\Omega\vert$, we still have $\Phi(\mu)=0$ but we do not have anymore $\Psi(\mu)=0$ and the property $(\ref{stab})$ is not satisfied. The state $\vert\Omega\rangle\langle\Omega\vert$ is no more an equilibrium state. As a consequence it is possible to observe more than one photon in the photon counter.\\

In the next section we deal with general strategy and the
particular problem of optimal control. Considerations about
optimal control is an interesting mean to point out the
importance of Markovian strategy.

\subsection{Optimal Control}

This section is then devoted to what is called the ``optimal
control'' problem. It deals with finding a particular control
strategy which must satisfy optimization constraint.
   In this section, we give the classical mathematical description of such problem and investigate general
    results in the discrete and in the continuous model of controlled quantum trajectories. Let us begin with the discrete model.

\subsubsection{The Discrete Case}

We come back to the description of a discrete quantum trajectory
for a two-level system as a Markov chain.

Let $n$ be fixed, thanks to Theorem $2$, a discrete controlled
quantum trajectory $(\rho_k^\mathbf{u})$ is described as follows.
Let $\rho$ be any state, if $\rho^\mathbf{u}_k=\rho$ then
$\rho^\mathbf{u}_{k+1}$ takes one of the values:
$$\mathcal{H}^{\mathbf{u},k}_i(\rho)=\frac{L_{i0}(k/n,u_{k}(n))(\rho)L^\star_{i0}(k/n,u_{k}(n))}
{Tr[L_{i0}(k/n,u_{k}(n))(\rho)L^\star_{i0}(k/n,u_{k}(n))]}\,\,i=0,1$$
with probability,
\begin{eqnarray*}p^\mathbf{u}_{k+1}(\rho)&=&Tr[L_{00}(k/n,u_{k}(n))(\rho)L^\star_{00}(k/n,u_{k}(n))]\,\,\textrm{for}\,\,i=0\\
 q_{k+1}^\mathbf{u}(\rho)&=&Tr[L_{10}(k/n,u_{k}(n))(\rho)L^\star_{10}(k/n,u_{k}(n))]\,\,\textrm{for}\,\,i=1.
\end{eqnarray*}

With this previous description, the property of a strategy $(u_k)$
can be enlarged. We can namely consider more general strategies
such that for all $k$ the term $u_k$ depend on all $(\rho_i)$ for
$i\leq k$. We define $\mathcal{U}$ the set of all admissible
strategies which satisfy this condition. Let us stress that in
this situation, the discrete quantum trajectory is no more a
Markov chain because the strategy at time $k$ depends on all the
past of the strategy.

 With this remark concerning the definition of strategies, we can expose the general problem
of ``optimal control''. In this article, we only consider finite
horizon problem. It is described as follows.

 Let $N$ be a fixed
integer and let $c$ and $\phi$ be two measurable function, the
optimal control problem in finite horizon is to consider what is
called the ``optimal cost'':
\begin{equation}
\min_{\mathbf{u}\in\mathcal{U}}\mathbf{E}\left[\sum_{k=0}^{N-1}
c(k,\rho^\mathbf{u}_k,u_{k})+\phi(\rho_{N}^\mathbf{u})\right].
\end{equation}

If there is some strategy which realizes the minimum, this
strategy is called the ``optimal strategy''. Let us investigate
the classical result in stochastic control for the optimal
strategy in this case.

For this we define:
$$V^k(\rho)=\min_{\mathbf{u}\in\mathcal{U}}\mathbf{E}
\left[\sum_{j=k}^{N-1}c(k,\rho^\mathbf{u}_j,u_{j})+
\phi(\rho^\mathbf{u}_{N})\bigg{/}\rho^\mathbf{u}_n=\rho\right].$$

\textbf{Remark} The function $c$ and $\phi$ are determined by the
optimization constraint imposed by the experience. The equation
which appears in the following theorem is called the cost equation
and the function $c$ and $\phi$ are called cost function.

\begin{thm}Let $\mathcal{U}$ be a compact set and suppose that $c$ is a continuous function. The solution of:
\begin{equation}\label{HJBd}\left\{\begin{array}{ccl}V^k(\rho)&=&\min_{u\in\mathcal{U}}\{p^\mathbf{u}_{k+1}(\rho)\mathcal{H}^{\mathbf{u},k}_0(\rho)+
q^\mathbf{u}_{k+1}(\rho)\mathcal{H}^{\mathbf{u},k}_1(\rho)+c(k,\rho,u_k)\}\\
V^N(\rho)&=&\phi(\rho)\end{array}\right.\end{equation} give the
optimal cost:
$$V^k(\rho)=\min_{\mathbf{u}\in\mathcal{U}}\mathbf{E}\left[\sum_{j=k}^{N-1}c(k,\rho^\mathbf{u}_j,u_{j})
+\phi(\rho_{N}^\mathbf{u})\bigg{/}\rho_n=\rho\right].$$ The
optimal strategy is given by:
\begin{equation}\label{opt}
u^\star:\rho\rightarrow u^\star_k(\rho)\in
\arg\min_{\mathbf{u}\in\mathcal{U}}\{p^\mathbf{u}_{k+1}(\rho)\mathcal{H}^{\mathbf{u},k}_0(\rho)
+q^\mathbf{u}_{k+1}(\rho)\mathcal{H}^{\mathbf{u},k}_1(\rho)+c(k,\rho,u_k)\}.
\end{equation}
Furthermore this strategy is Markovian.
\end{thm}

\begin{pf}
The proof is based of what is called dynamic programming in
stochastic control theory. Let $\mathbf{u}$ be any strategy and
let $V$ defined by the formula $(\ref{HJBd})$, we
have\begin{eqnarray*}&&\mathbf{E}[(V^{k+1}\left(\rho^\mathbf{u}_{k+1}\right)-V^{k}\left(\rho^\mathbf{u}_k)\right)/\sigma\{\rho_i^\mathbf{u},i\leq
k\}]
\\&=&p_{k+1}^\mathbf{u}V^{k+1}\left(\mathcal{H}^{\mathbf{u},k}_0(\rho_k^\mathbf{u})\right)+q_{k+1}^\mathbf{u}V^{k+1}\left(\mathcal{H}^{\mathbf{u},k}_1(\rho_k^\mathbf{u})\right)
-V^{k}(\rho^\mathbf{u}_k))
\end{eqnarray*}
 then we have
\begin{eqnarray*}
 &&\mathbf{E}\left[V^N(\rho^\mathbf{u}_N)-V^0(\rho)\right]\\
&=&\sum_{k=0}^{N-1}\mathbf{E}\left[V^{k+1}(\rho^\mathbf{u}_{k+1})-V^k(\rho^\mathbf{u}_k)\right]\\
&=&\sum_{k=0}^{N-1}\mathbf{E}\left[p_{k+1}^\mathbf{u}V^{k+1}\left(\mathcal{H}^{\mathbf{u},k}_0(\rho_k^\mathbf{u})\right)+q_{k+1}^\mathbf{u}V^{k+1}\left(\mathcal{H}^{\mathbf{u},k}_1(\rho_k^\mathbf{u})\right)
-V^{k}(\rho^\mathbf{u}_k))\right]\\
&\geq&-\sum_{k=0}^{N-1}\mathbf{E}\left[c(k,\rho_k^\mathbf{u},u_k)\right]\,\,\,(\textrm{by
definition of the min}).
\end{eqnarray*}
Hence for all strategy $\mathbf{u}$, we have
$$V^0(\rho)\leq\mathbf{E}\left[\sum_{k=0}^{N-1}c(k,\rho_k^\mathbf{u},u_k)+\phi(\rho_N^\mathbf{u})\right].$$

Moreover we have equality if we choose the strategy given by the
formula $(\ref{opt})$.
 This strategy is Markovian because the function $c$ depends only on $\rho_k$ at time $k$.
\end{pf}\\

The system $(\ref{HJBd})$ which describes the cost equation is
called the discrete Hamilton-Jacobi Bellman equation.

The fact that the optimal strategy is Markovian is another
justification of the choice of such model of control for the
discrete quantum trajectory. This theorem claims that we need just
Markovian strategy in order to solve the ``optimal control''
problem.

The next last section is devoted to the same investigation in
continuous time models of quantum trajectories.

\subsubsection{The Continuous Case}

In the third section, we have proved the Poisson and the diffusion
approximation in quantum measurement theory. We have the diffusive
evolution equation
\begin{eqnarray}\label{diffff}
 \rho_t&=&\rho_0+\int_0^tL(s,\rho^\mathbf{u}_s,u(s,\rho^\mathbf{u}_s))ds+
 \int_{0}^t\Theta(s,\rho_s^\mathbf{u},u(s,\rho^\mathbf{u}_s))dW_s,
\end{eqnarray}
and the jump-equation is
\begin{eqnarray}\label{poissss}
 \rho_t&=&\rho_0+\int_0^tR(s\mbox{\tiny{$-$}},\rho^\mathbf{u}_{s\mbox{\tiny{$-$}}},u(s\mbox{\tiny{$-$}},\rho^\mathbf{u}_{s\mbox{\tiny{$-$}}}))ds
 \\&&+\int_0^t\int_{\mathbb{R}}Q(s,\rho^\mathbf{u}_{s\mbox{\tiny{$-$}}},u(s\mbox{\tiny{$-$}},\rho^\mathbf{u}_{s\mbox{\tiny{$-$}}}))
 \mathbf{1}_{0<x<Tr[\mathcal{J}
 (s\mbox{\tiny{$-$}},u(s\mbox{\tiny{$-$}}
 ,\rho^\mathbf{u}_{s\mbox{\tiny{$-$}}}))(\rho^\mathbf{u}_{s\mbox{\tiny{$-$}}})]}N(dx,ds),
 \end{eqnarray}
where the functions $L$, $\Theta$, $R$ and $Q$ are defined in
Section $2$.

In this section, we consider the same problem of "optimal control"
as in the discrete case. Let
$(\Omega,\mathcal{F},\mathcal{F}_t,P)$ be a probability space
where we consider the diffusive equation
\begin{eqnarray*}
 \rho_t&=&\rho_0+\int_0^tL(s,\rho^\mathbf{u}_s,u_s)ds+\int_{0}^t\Theta(s,\rho_s^\mathbf{u},u_s)dW_s,
 \end{eqnarray*}
 and the jump-equation
 \begin{eqnarray*}
 \rho_t&=&\rho_0+\int_0^tR(s\mbox{\tiny{$-$}},\rho^\mathbf{u}_{s\mbox{\tiny{$-$}}},u_{s\mbox{\tiny{$-$}}})ds\\&&+\int_0^t\int_{\mathbb{R}}Q(s,\rho^\mathbf{u}_{s\mbox{\tiny{$-$}}},u_{s\mbox{\tiny{$-$}}})
 \mathbf{1}_{0<x<Tr[\mathcal{J}(s\mbox{\tiny{$-$}},u_{s\mbox{\tiny{$-$}}})(\rho^\mathbf{u}_{s\mbox{\tiny{$-$}}})]}N(dx,ds)
 \end{eqnarray*}
 where the strategy $\mathbf{u}=(u_t)$ is just supposed to be a function $\mathcal{F}_t$ adapted (not only Markovian). In the case where $\mathcal{F}_t$ corresponds to the filtration generated by the process $(\rho_t)$, we recover the same definition as the discrete case.
  Concerning existence and uniqueness of a solution, with the condition $(\ref{localcondition})$ of Section $2.1$ for the functions $L$, $R$ and $\theta$ the previous
  equations
admit a unique solution. Furthermore the solution takes values in
the set of states on $\mathcal{H}_0$. The set of all admissible
strategy which satisfy the condition of adaptation is also denoted
by $\mathcal{U}$. The optimal control problem in this situation is
expressed as follows.

Let $c$ and $\phi$ be two cost function. Let $T>0$, the optimal
problem in finite horizon is given by
\begin{equation}
 \min_{\mathbf{u}\in\mathcal{U}}\mathbf{E}\left[\int_0^Tc(s,\rho_s^\mathbf{u},u_s)ds+\phi(\rho_T^\mathbf{u})\right].
\end{equation}
 As in the discrete model, we introduce the following function:
\begin{equation}\label{coo}
 V(t,\rho)=\min_{\mathbf{u}\in\mathcal{U}}\mathbf{E}
 \left[\int_t^Tc(s,\rho_s^\mathbf{u},u_s)ds+\phi
 (\rho_T^\mathbf{u})\bigg{/}\rho^\mathbf{u}_t=\rho\right],
\end{equation}
which satisfies
$$V(0,\rho_0)=\min_{\mathbf{u}\in\mathcal{U}}\mathbf{E}\left[\int_0^Tc(s,\rho_s^\mathbf{u},u_s)ds+\phi(\rho_T^\mathbf{u})\right].$$
 The function $(\ref{coo})$ represents the result of optimal
control after $t$ assuming $\rho_t=\rho$.

In this article, we just give the result for the optimal control
problem for the diffusive case. A similar result for the Poisson
case can be found in \cite{bouten-2005}.

As in the discrete case, it appears a continuous time version of
the Hamilton-Jacobi-Bellmann Equation. The usual expression of
this equation use the notion of infinitesimal generator of
$(\rho^\mathbf{u}_t)$. It is described as follows in our context.
 A quantum trajectory $(\rho^\mathbf{u}_t)$ is considered as a process which takes values in $\mathbb{R}^3$
  with the identification of the state and the Bloch sphere $\mathbb{B}_1(\mathbb{R}^3)=\{(x,y,z)\in\mathbb{R}^3/x^2+y^2+z^2\leq1\}$, that
  is,
$$\begin{array}{cccc}
\Phi:&\mathbb{B}_1(\mathbb{R}^3)&\longmapsto&\mathbb{M}_2(\mathbb{C})\\
&(x,y,z)&\longrightarrow&\frac{1}{2}\left(\begin{array}{cc}1+x&y+iz\\y-iz&1-x\end{array}\right)
  \end{array}$$
 The map $\Phi$ is injective and its range is the set of states. By considering that the functions $L$ and $\Theta$ are
 applications
 from $\mathbb{R}_+\times\mathbb{R}^3$ to $\mathbb{R}^3$,
  the stochastic differential equation concerning the diffusive case can be written
   as a system of stochastic differential equation on $\mathbb{R}^3$ of the form:
$$(\rho^\mathbf{u}_t)_i=\rho_0+\int_0^tL_i(s,\rho^\mathbf{u}_s,u_s)ds+\int_{0}^t\Theta_i(s,\rho_s^\mathbf{u},u_s)dW_s\,\,\,\,i\in\{1,2,3\}$$
 where $(\rho^\mathbf{u}_t)_i$
(respectively $\Theta_i$ and $L_i$) corresponds to the coordinate
function of $\rho^\mathbf{u}_t$ (respectively $\Theta$ and $L$).

 We introduce the $3\times 3$
matrix $\Pi$ defined by $\Pi_{ij}=\Theta_i\Theta_j$. The
infinitesimal generator $\mathcal{A}^{u,t}$ of the process
$(\rho^\mathbf{u}_t)$ acts on the functions $f$ which are $C^2$
and bounded in the following way
\begin{equation}\label{gene}
 \mathcal{A}^{u,t}f(x)=\frac{1}{2}\sum_{i,j=1}^3\Pi_{ij}(t,x,u)\frac{\partial f(x)}{\partial x_i\partial x_j}+\sum_{i=1}^3L_i(t,x,u)\frac{\partial f(x)}{\partial x_i}.
\end{equation}
for all $t\geq0$, $u\in\mathbb{R}$ and $x\in\mathbb{R}^3$. In
particular if $u$ is a fixed constant, let $(\rho_t)$ be the
solution of
$$(\rho_t)_i=\rho_0+\int_0^tL_i(s,\rho_s,u)ds+\int_{0}^t\Theta_i(s,\rho_s,u)dW_s\,\,\,\,i\in\{1,2,3\}.$$ Hence for all function $f$ which is $C^2$ and bounded,
the following process
$$\mathcal{M}_t=f(\rho_t)-f(\rho_0)-\int_{0}^t\mathcal{A}^{u,s}f(\rho_s)ds$$
 is a martingale for the filtration generated by $(\rho_t)$.

The following theorem express the result in optimal control for
the diffusive quantum trajectory.

\begin{thm}
 Suppose there is a function $(t,\rho)\rightarrow V(t,\rho)$ which is $C^1$ in $t$ and $C^2$ in $\rho$ such that:
\begin{equation}\label{HJBCONT}
 \left\{\begin{array}{ccl}\frac{\partial V(t,\rho)}{\partial t}+\min_{u\in\mathcal{U}}\{\mathcal{A}^{u,t}V(t,\rho)+c(t,\rho,u)\}&=&0\\
V(T,\rho)&=&\phi(\rho)\end{array}\right.
\end{equation}
where $\mathcal{A}^{u,t}f(x)$ is defined by the expression
$(\ref{gene})$. The function $V$ gives the solution of the optimal
problem, that is,
$$V(t,\rho)=\min_{\mathbf{u}\in\mathcal{U}}\mathbf{E}\left[\int_t^Tc(s,\rho_s^\mathbf{u},u_s)ds+
\phi(\rho_T^\mathbf{u})\bigg{/}\rho^\mathbf{u}_t=\rho\right].$$
Furthermore if the strategy $\mathbf{u}$ defined by
\begin{equation}
u^\star(t,\rho)\,\in\,\arg
\min_{u\in\mathcal{U}}\{\mathcal{A}^{u,t}V(t,\rho)+c(t,\rho,u)\}
\end{equation}
is an admissible strategy then it defines an optimal strategy.
Moreover this strategy is Markovian.
\end{thm}

The equation $(\ref{HJBCONT})$ is the Hamilton-Jacobi-Bellmann
equation in the continuous case.

A proof of this theorem can be found in \cite{MR0280248} or
\cite{pham-2005-2}. The interest of such theorem in our context is
to notice that the optimal strategy is Markovian, this confirms
the choice of such strategy in the model of quantum trajectories
with control.

A similar result holds for the Poisson case. The infinitesimal
generator for such process is given in \cite{cheridito-2005-15},
explicit computations can be found in \cite{multi}.

\end{document}